\documentclass[12pt]{article}

\usepackage{amssymb}
\usepackage{amsmath}
\usepackage{cite}
\usepackage[arrow, matrix, curve]{xy}
\usepackage{bbm}

\usepackage[dvips]{graphicx}
\usepackage{xcolor}

\begin{document}

\renewcommand{\citeleft}{{\rm [}}
\renewcommand{\citeright}{{\rm ]}}
\renewcommand{\citepunct}{{\rm,\ }}
\renewcommand{\citemid}{{\rm,\ }}

\newcounter{abschnitt}
\newtheorem{satz}{Theorem}
\newtheorem{coro}[satz]{Corollary}
\newtheorem{theorem}{Theorem}[abschnitt]
\newtheorem{koro}[theorem]{Corollary}
\newtheorem{prop}[theorem]{Proposition}
\newtheorem{lem}[theorem]{Lemma}
\newtheorem{expls}[theorem]{Examples}

\newcommand{\mres}{\mathbin{\vrule height 1.6ex depth 0pt width 0.11ex \vrule height 0.11ex depth 0pt width 1ex}}

\renewenvironment{quote}{\list{}{\leftmargin=0.62in\rightmargin=0.62in}\item[]}{\endlist}

\newcounter{saveeqn}
\newcommand{\alpheqn}{\setcounter{saveeqn}{\value{abschnitt}}
\renewcommand{\theequation}{\mbox{\arabic{saveeqn}.\arabic{equation}}}}
\newcommand{\reseteqn}{\setcounter{equation}{0}
\renewcommand{\theequation}{\arabic{equation}}}

\hyphenpenalty=9000

\sloppy

\phantom{a}

\vspace{-1.7cm}

\begin{center}
\begin{Large} {\bf Fixed Points of Minkowski Valuations} \\[0.4cm] \end{Large}

\begin{large} Oscar Ortega-Moreno and Franz E.\ Schuster\end{large}
\end{center}

\vspace{-0.9cm}

\begin{quote}
\footnotesize{ \vskip 1cm \noindent {\bf Abstract.}
It is shown that for any sufficiently regular even Minkowski valuation $\Phi$ which is homogeneous and intertwines rigid motions, there exists a neighborhood of the unit ball, where balls are the only solutions to the fixed-point problem $\Phi^2 K = \alpha K$. This significantly generalizes results by Ivaki for projection bodies and suggests, via the Lutwak--Schneider class reduction technique, a new approach to Petty's conjectured projection inequality.
}
\end{quote}

\vspace{0.4cm}

\centerline{\large{\bf{ \setcounter{abschnitt}{1}
\arabic{abschnitt}. Introduction}}}

\alpheqn

\vspace{0.5cm}

The classical Brunn--Minkowski theory has emerged from Minkowski's studies of the volume of sums of convex bodies, leading to such central notions as intrinsic volumes or, more generally, mixed volumes and the basic inequalities between them. It was also Minkowski who introduced projection bodies of convex bodies -- another core concept of the theory -- which were later discovered to be objects of independent investigations in several areas (see, e.g., the excellent monographs by Gardner \textbf{\cite{gardner2ed}}, Koldobsky \textbf{\cite{KoldobskyBook}}, and Schneider \textbf{\cite{schneider93}}).
For their precise definition, let $K$ be a convex body (that is, a compact, convex set) in $\mathbb{R}^n$, where throughout $n \geq 3$, and recall that $K$ is determined uniquely by its support function $h(K,u) = \max \{u \cdot x: x \in K\}$, $u \in  \mathbb{S}^{n-1}$. The \emph{projection body} $\Pi K$ of $K$ is the convex body defined by
\[h(\Pi K,u) = V_{n-1}(K|u^{\perp}), \qquad u \in \mathbb{S}^{n-1}.  \]
Here, $V_i$, $0 \leq i \leq n$, denotes the $i$-th intrinsic volume (see Section~2 for definition) and $K|u^{\perp}$ is the orthogonal projection onto the hyperplane $u^{\perp}$.

The polar projection inequality of Petty \textbf{\cite{Petty1972}}, providing the classical relation between the volume of a convex body and its \emph{polar} projection body, is an affine invariant inequality that not only significantly improves the classical isoperimetric inequality but had a tremendous impact in geometric analysis that can still be felt to this day (see, e.g., \textbf{\cite{boeroezky2013, habschu09, HaberlSchu2019, LYZ2000a, LYZ2010a}} and the references therein). In view of the Blaschke--Santal\'o inequality (see, e.g., \textbf{\cite[\textnormal{Section 10.7}]{schneider93}}), an analog of Petty's inequality for the volume of projection bodies (as opposed to that of polar projection bodies) would provide an even stronger affine isoperimetric inequality. However, this remains one of the major open problems in convex geometric analysis.

\vspace{0.25cm}

\noindent {\bf Conjecture} (Petty \textbf{\cite{Petty1972}}) \emph{Among convex bodies $K \subseteq \mathbb{R}^n$ of non-empty interior, the volume ratio $V_n(\Pi K)/V_n(K)^{n-1}$ is minimized precisely by ellipsoids.}

\vspace{0.25cm}

Despite the universally acknowledged importance of Petty's conjecture, progress has been slow over the last decades
(see, e.g., \textbf{\cite{GiannopolousPapa, Ivaki2018, Lutwak1990a, Saroglou2015, SaroglouZvavitch, Schneider1987}}). Especially indicative of the difficulty of the problem is a recent disheartening observation by Saroglou \textbf{\cite{Saroglou2011}} that in general, the volume of projection bodies does not decrease under the standard symmetrization technique -- Steiner symmetrization -- used to establish many powerful isoperimetric inequalities in geometry.

\pagebreak

As Saroglou's observation demonstrates, Petty's conjecture requires an approach different from Steiner symmetrization. Before discussing one such possibility, we want to mention here that Lutwak showed that if Petty's conjecture is correct, then it would imply a whole family of isoperimetric inequalities for projection bodies of different degrees. To make this more precise, recall that for $1 \leq i \leq n - 1$, the projection body of degree $i$, $\Pi_iK$, of a convex body $K \subseteq \mathbb{R}^n$ is defined by
\[h(\Pi_i K,u) = V_{i}(K|u^{\perp}), \qquad u \in \mathbb{S}^{n-1}.  \]
Note that, clearly, $\Pi_{n-1} = \Pi$ is just Minkowski's projection body map.

In 1993, Lutwak \textbf{\cite{Lutwak1993}} showed how the polar Petty projection inequality can be used to deduce analogous inequalities for polar projection bodies of all degrees. Previously, Lutwak \textbf{\cite{Lutwak1990b}} had conjectured that among bodies $K \subseteq \mathbb{R}^n$ of non-empty interior, the ratios $V_{i+1}(\Pi_i K)/V_{i+1}(K)^i$ are minimized (precisely) by Euclidean balls for all $1 \leq i \leq n - 1$ and had shown that this would follow if Petty's conjecture holds true. Moreover, the Lutwak--Petty conjecture was confirmed in  \textbf{\cite{Lutwak1990b}} when $i = 1$ by applying the following instance of the \emph{class reduction technique}, first observed in the context of Petty's conjecture by Schneider and later extended by Lutwak.

\begin{prop} \label{LutwakSchnProp} \textnormal{(\!\!\textbf{\cite{Lutwak1990b}}, \textbf{\cite{Schneider1987}})}
If $1 \leq i \leq n - 1$ and $K \subseteq \mathbb{R}^n$ is a convex body whose dimension is at least $i + 1$, then
\[ \frac{V_{i+1}(\Pi_i K)}{V_{i+1}(K)^i} \geq \frac{V_{i+1}(\Pi_i^2 K)}{V_{i+1}(\Pi_i K)^i} \]
with equality if and only if $\Pi_i^2 K$ and $K$ are homothetic.
\end{prop}

Note that Proposition \ref{LutwakSchnProp} implies that if minimizers of $V_{i+1}(\Pi_i K)/V_{i+1}(K)^i$ do exist, then they must be translates of solutions to the fixed-point problem
\begin{equation} \label{fixedpointPii}
\Pi_i^2K = \alpha K
\end{equation}
for some $\alpha > 0$. When $i = n - 1$, minimizers of $V_n(\Pi K)/V_n(K)^{n-1}$ do exist due to affine invariance and ellipsoids are solutions to (\ref{fixedpointPii}). However, it is well known that the projection body map $\Pi$ also admits a large class of polytopal solutions to (\ref{fixedpointPii}) which were completely classified by Weil \textbf{\cite{Weil1971}}. This is in stark contrast to the case $i = 1$, where Schneider \textbf{\cite{Schneider1977}} showed that $\Pi_1^2K$ is homothetic to $K$ if and only if $K$ is a ball. Thus, after establishing the existence of minimizers for $V_2(\Pi_1 K)/V_2(K)$, Lutwak \textbf{\cite{Lutwak1990b}} concluded that they are precisely the Euclidean balls.

While the fixed-point problem (\ref{fixedpointPii}) can be   resolved easily for $i = 1$, it is a much harder problem for $i > 1$ (cf.\ Section 6). However, a breakthrough was achieved by Ivaki \textbf{\cite{Ivaki2017, Ivaki2018}} recently, who proved that there exists a $C^2$ neighborhood of the unit ball $\mathbb{B}^n$, where, for $2 \leq i \leq n - 2$, Euclidean balls  and, for $i = n - 1$, ellipsoids are the only solutions to (\ref{fixedpointPii}). The case $i = n - 1$ of Ivaki's results can also be deduced from an independent stronger theorem of Saroglou and Zvavitch \textbf{\cite{SaroglouZvavitch}} about iterations of $\Pi$ which also confirms Petty's conjecture in a $C^2$~neighborhood of the unit ball.

\pagebreak

The main goal of this article is to generalize Ivaki's results on local solutions of (\ref{fixedpointPii}) to a large class of Minkowski valuations. A \emph{Minkowski valuation} on the space $\mathcal{K}^n$ of convex bodies in $\mathbb{R}^n$ endowed with the Hausdorff metric is a map $\Phi: \mathcal{K}^n \rightarrow \mathcal{K}^n$ such that
\[\Phi(K) + \Phi(L) = \Phi(K \cup L) + \Phi(K \cap L)  \]
whenever $K \cup L$ is convex and addition on $\mathcal{K}^n$ is the usual vector addition. Scalar valued valuations have long been at the center of convex and integral geometry
(see, e.g., \textbf{\cite{Alesker1999, Alesker2001, Aleskerbook2018, BernigFu2011, HabParap2014, Klain:Rota, LudwigReitz2010}}). Their systematic study was initiated by Hadwiger~\textbf{\cite{hadwiger51}}, whose characterization of the intrinsic volumes as the continuous rigid motion invariant scalar valuations is one of the most famous results in valuation theory.

First results on Minkowski valuations were obtained in 1974 by Schneider~\textbf{\cite{Schneider1974c}}. However, they became the focus of increased attention only after the seminal work of Ludwig \textbf{\cite{ludwig02, Ludwig:Minkowski}} on Minkowski valuations intertwining affine transformations. For example, Ludwig established a characterization of Minkowski's projection body map as the unique continuous translation invariant Minkowski valuation which is contravariant with respect to volume preserving linear transformations.

The contributions of several authors (see, e.g., \textbf{\cite{abardber11, boeroezkyludwig2019, colesantietal2017, haberl11, ludwig2010, SchuWan12, Wannerer2011}}) show that Minkowski valuations compatible with affine  transformations often form convex cones generated by \emph{finitely} many maps. In contrast, Minkowski valuations intertwining rigid motions form an \emph{infinite dimensional} cone, containing, e.g., the projection bodies of arbitrary degrees. The efforts to obtain an analogue of Hadwiger's theorem for Minkowski valuations
(see \textbf{\cite{kiderlen05, Schu06a, Schu10, SchuWan13, SchuWan16}}) culminated in the recent work of Dorrek \textbf{\cite{Dorrek2017b}} who established the following \emph{spherical convolution} (see Section 2) representation under the additional assumption of homogeneity. Throughout, a map $\Phi: \mathcal{K}^n \rightarrow \mathcal{K}^n$ is said to have degree $i$ if $\Phi(\lambda K) = \lambda^i \Phi(K)$ for all $K \in \mathcal{K}$ and $\lambda \geq 0$. (By a result of McMullen, any translation invariant continuous valuation that is also homogeneous must be of integer degree $i \in \{0, \ldots, n\}$.)

\begin{theorem} \label{thmdorrek} \textnormal{(\!\!\!\textbf{\cite{Dorrek2017b}})}
If $\Phi_i: \mathcal{K}^n \rightarrow \mathcal{K}^n$ is a continuous translation invariant Minkowski valuation of degree $1 \leq i \leq n - 1$ which commutes with $\mathrm{SO}(n)$, then there exists a unique
 $\mathrm{SO}(n - 1)$ invariant $f \in L^1(\mathbb{S}^{n-1})$ with center of mass at the origin such that for every $K \in \mathcal{K}^n$,
\begin{equation} \label{convrep}
h(\Phi_i K,\cdot) = S_i(K,\cdot) \ast f.
\end{equation}
\end{theorem}

\vspace{0.2cm}

The measures $S_i(K,\cdot)$, $1 \leq i \leq n - 1$, on $\mathbb{S}^{n-1}$ are the \emph{area measures of order~$i$} associated with $K$ (see Section 2 for definition).
We call the function $f$ in (\ref{convrep}) the \emph{generating function} of $\Phi_i$. While a complete classification of all such generating functions has not yet been obtained, it is known that
for any $i$, we may take $f$ in (\ref{convrep}) to be the support function of an arbitrary body of revolution $L \in \mathcal{K}^n$. In this case, we say that $\Phi_i$ is generated by $L$.
If, in addition, the boundary of $L$ is a $C^2$ submanifold of $\mathbb{R}^n$ with everywhere positive Gaussian curvature, we call $\Phi_i$ a \emph{$C^2_+$~regular} Minkowski valuation.
If $\Phi_i K = \{o\}$ for all $K \in \mathcal{K}^n$, we call $\Phi_i$ \emph{trivial}.

\pagebreak

Over the past 15 years, it has become more and more apparent that several classic inequalities involving projection bodies (of arbitrary degree) hold, in fact, for the entire class, or at least a large subclass, of Minkowski valuations intertwining rigid motions (see, e.g., \textbf{\cite{ABS2011, BPSW2014, HaberlSchu2019, HofstaetterSchu2021, papschu12, Schu06a}}). Among the first results in this
direction, it was proved in \textbf{\cite{Schu06a}} that if $\Phi_1: \mathcal{K}^n \rightarrow \mathcal{K}^n$ is a non-trivial continuous translation
invariant Minkowski valuation of degree $1$ which commutes with $\mathrm{SO}(n)$ and is \emph{monotone} w.r.t.\ set inclusion, then
$V_{2}(\Phi_1 K)/V_{2}(K)$ is minimized, among $K \in \mathcal{K}^n$ with non-empty interior, precisely by Euclidean balls. This generalized not only Lutwak's inequality for
$\Pi_1$ but also strongly relied on the class reduction technique which, in this case, was shown in \textbf{\cite{Schu06a}} to imply that the minimizers must be translates of the solutions to the fixed-point problem $\Phi_1^2 K = \alpha K$ for some $\alpha > 0$. The latter had previously been solved under the assumption of monotonicity by
Kiderlen \textbf{\cite{kiderlen05}}, who showed that balls are the only solutions.

Our first goal is to establish the following extension of the Lutwak--Petty class reduction from Proposition \ref{LutwakSchnProp} to all monotone homogeneous Minkowski valuations intertwining
rigid motions (for the degree 1 case, see also \textbf{\cite{Schu06a}}).

\vspace{0.3cm}

\noindent {\bf Proposition 1} \emph{Let $1 \leq i \leq n - 1$ and $\Phi_i: \mathcal{K}^n \rightarrow \mathcal{K}^n$ be a non-trivial monotone and translation invariant Minkowski valuation of degree $i$ which commutes with $\mathrm{SO}(n)$. If $K \in \mathcal{K}^n$ has dimension at least $i + 1$, then
\begin{equation} \label{prop1inequ}
\frac{V_{i+1}(\Phi_i K)}{V_{i+1}(K)^i} \geq \frac{V_{i+1}(\Phi_i^2 K)}{V_{i+1}(\Phi_iK)^i}
\end{equation}
with equality if and only if $\Phi_i^2 K$ and $K$ are homothetic. Moreover, if $\Phi_i$ is $C_2^+$~regular and balls are the only solutions to the fixed-point problem $\Phi_i^2 K = \alpha K$ for some $\alpha > 0$, then $V_{i+1}(\Phi_i K)/V_{i+1}(K)^i$ is minimized precisely by Euclidean balls.}

\vspace{0.3cm}

Let us emphasize how Proposition 1 suggests a new approach towards Petty's conjecture that has the advantage of introducing higher regularity and therefore, in particular, eliminates possible polytopal solutions of the associated fixed-point problems. In order to describe this approach, first note that in \textbf{\cite{Schu06a}} it was shown that every continuous
Minkowski valuation $\Phi_{n-1}$ intertwining rigid motions of degree $n - 1$ which is \emph{even} (that is, $\Phi_{n-1} K = \Phi_{n-1}(-K)$ for all $K \in \mathcal{K}^n$) is generated by a body of revolution $L \in \mathcal{K}^n$ and, thus, is monotone (cf.\ Section 3). Consequently, by Proposition~1, choosing suitable even Minkowski valuations $\Phi_{n-1}$ generated by (sufficiently) smooth bodies $L \in \mathcal{K}^n$ approximating the segment -- the generating body of $\Pi$ -- and showing that the only solutions to $\Phi_{n-1}^{2}K = \alpha K$ are Euclidean balls, would confirm Petty's conjecture (at least, up to equality conditions). The fundamental difference to the original class reduction approach is that for Minkowski valuations $\Phi_{n-1}$ generated by bodies $L$ of, say, class $C^2_+$, $\Phi_{n-1}K$ and, therefore, $\Phi_{n-1}^2K$ will also belong to the class $C^2_+$ for \emph{all} $K \in \mathcal{K}^n$. In particular, the possible minimizers of $V_n(\Phi_{n-1}K)/V_n(K)^{n-1}$ must also be of class $C^2_+$.

With our main theorem, we generalize Ivaki's results \textbf{\cite{Ivaki2017, Ivaki2018}} about solutions to the fixed-point problems $\Pi_i^2 K = \alpha K$ locally around balls to all sufficiently regular even Minkowski valuations intertwining rigid motions of a fixed degree $i \in \{2, \ldots, n - 1\}$ (the case $i = 1$ having been \emph{globally} settled by Kiderlen \textbf{\cite{kiderlen05}}).

\begin{satz} \label{mainthm} Let $2 \leq i \leq n - 1$ and $\Phi_i: \mathcal{K}^n \rightarrow \mathcal{K}^n$ be a $C^2_+$ regular translation invariant even Minkowski valuation of degree $i$ which commutes with $\mathrm{SO}(n)$.
Then there exists $\varepsilon > 0$ such that if $K \in \mathcal{K}^n$ has a $C^2$ support function and satisfies
\begin{itemize}
\item[(i)] $\|h(\gamma K + x,\cdot) - h(\mathbb{B}^n,\cdot)\|_{C^2(\mathbb{S}^{n-1})} < \varepsilon$ for some $\gamma > 0$ and $x \in \mathbb{R}^n$,
\item[(ii)] $\Phi_i^2 K = \alpha K$ for some $\alpha > 0$,
\end{itemize}
then $K$ must be a Euclidean ball.
\end{satz}

In Section 6 we will in fact prove a more general result, Theorem \ref{mostgeneralthm}, than Theorem \ref{mainthm}. It provides sufficient conditions on the generating function
$f \in  L^1(\mathbb{S}^{n-1})$ of a general continuous even Minkowski valuation $\Phi_i$ intertwining rigid motions, to conclude that in a neighborhood around the ball, the only solutions
to the fixed-point problem $\Phi_i^2 K = \alpha K$ are Euclidean balls. These conditions are stated in terms of the spherical harmonic expansion of $f$ and are easily checked
to be satisfied by the support function of the segment. In this way, we generalize Ivaki's results \textbf{\cite{Ivaki2017}} for projection bodies of order $2 \leq i \leq n - 2$. The case
$i = n - 1$ is not covered due to affine contravariance of $\Pi$, but Ivaki's proof \textbf{\cite{Ivaki2018}} relies on similar techniques.

In order to prove Theorem \ref{mainthm}, we confirm that the conditions on generating functions contained in Theorem \ref{mostgeneralthm} are satisfied by support functions of
convex bodies of revolution of class $C^2_+$. This boils down to a novel spectral gap theorem of independent interest for spherical convolution operators generated
by such support functions which we obtain in Section 5. Moreover, we also discuss in Section 6 that when $n \geq 4$ and $2 \leq i \leq n - 2$, the conditions from Theorem \ref{mostgeneralthm} are satisfied for even Minkowski valuations generated by generalized zonoids of revolution.

\vspace{1cm}

\centerline{\large{\bf{ \setcounter{abschnitt}{2}
\arabic{abschnitt}. Background material}}}

\reseteqn \alpheqn \setcounter{theorem}{0}

\vspace{0.6cm}

In the following we first recall basic facts about convex bodies, mixed volumes, and area measures. Next, we collect the required material from
harmonic analysis, in particular, about the spherical convolution of measures and its relation to the theory of spherical harmonics. In the final part of this section, we recall
the definition of Frech\'et derivatives and state a useful version of the inverse function theorem.
General references for the material of this section are the monographs by Gardner \textbf{\cite{gardner2ed}}, Schneider \textbf{\cite{schneider93}}, and Groemer \textbf{\cite{Groemer1996}},
as well as the celebrated exposition \textbf{\cite{Hamilton1982}} on the inverse function theorem by Hamilton.

First recall that each $K \in \mathcal{K}^n$ is uniquely determined by its support function $h(K,x) = \max \{x \cdot y: y \in K\}$, $x \in \mathbb{R}^n$, which is (positively) homogeneous of degree one and subadditive. Conversely, every function on $\mathbb{R}^n$ with these two properties is the support function of a unique body in $\mathcal{K}^n$. In particular, a function $h \in C^2(\mathbb{R}^n)$ which is homogeneous of degree one is the support function of a convex body $K \in \mathcal{K}^n$ if and only if its Hessian $D^2h(u)$ is positive semi-definite for all $u \in \mathbb{S}^{n-1}$.

For $K, L \in \mathcal{K}^n$, their Minkowski sum is given by $K + L = \{x + y: x \in K, y \in L\}$ and its support function by $h(K + L,\cdot) = h(K,\cdot) + h(L,\cdot)$.
Moreover, for every $\vartheta \in \mathrm{SO}(n)$ and $y \in \mathbb{R}^n$, we have
\begin{equation} \label{sontranslsuppfct}
h(\vartheta K,x) = h(K,\vartheta^{-1}x) \qquad \mbox{and} \qquad h(K + y,x) = h(K,x) + x \cdot y
\end{equation}
for all $x \in \mathbb{R}^n$. Next, recall that the Hausdorff distance $d(K,L)$ of $K, L \in \mathcal{K}^n$ can be expressed by $d(K,L) = \|h(K,\cdot) - h(L,\cdot)\|_{\infty}$,
where $\|\,\cdot\,\|_{\infty}$ denotes the maximum norm on $C(\mathbb{S}^{n-1})$. Moreover, $K \subseteq L$ if and only if $h(K,\cdot) \leq h(L,\cdot)$, in particular, $h(K,\cdot) > 0$ if and only if $o \in \mathrm{int}\,K$.

A body $K \in \mathcal{K}^n$ is said to be of class $C^k_+$ if its boundary hypersurface $\partial K$ is a $C^k$ submanifold of $\mathbb{R}^n$ and the map $n_K:\partial K \to \mathbb{S}^{n-1}$ that maps a boundary point to its unique outer unit normal is a $C^k$ diffeomorphism. Equivalently, $K \in \mathcal{K}^n$ is of class $C^k_+$ if $h(K,\cdot) \in C^k(\mathbb{R}^n)$ and the restriction of the Hessian $D^2h(K,\cdot)(u)$ to $u^{\perp}$ is positive
definite for every $u \in \mathbb{S}^{n-1}$.

By a classical result of Minkowski, the volume of a Minkowski linear combination $\lambda_1K_1 + \cdots + \lambda_mK_m$, where $K_1, \ldots, K_m \in \mathcal{K}^n$ and $\lambda_1, \ldots,
\lambda_m \geq 0$, can be expressed as a homogeneous polynomial of degree $n$,
\begin{equation} \label{mixed}
V_n(\lambda_1K_1 + \cdots +\lambda_m K_m)=\sum \limits_{j_1,\ldots, j_n=1}^m V(K_{j_1},\ldots,K_{j_n})\lambda_{j_1}\cdots\lambda_{j_n},
\end{equation}
where the coefficients $V(K_{j_1},\ldots,K_{j_n})$ are the {\it mixed volumes} of $K_{j_1}, \ldots, K_{j_n}$ which depend only on $K_{j_1}, \ldots, K_{j_n}$ and are symmetric in their
arguments. Moreover, mixed volumes are translation invariant, Minkowski additive, monotone w.r.t.\ set inclusion in each of their arguments, and $V(K_1,\ldots,K_n) > 0$ if and only if there are segments
$l_i \subseteq K_i$, $1 \leq i \leq n$, with linearly independent directions.

For $K, L \in \mathcal{K}^n$ and $0 \leq i \leq n$, let $V(K[i],L[n-i])$ denote the mixed volume with $i$ copies of $K$ and $n - i$ copies of $L$. The \emph{$i$th intrinsic volume} of $K$ is given by
\begin{equation*} \label{viwi}
V_i(K)=\frac{1}{\kappa_{n-i}}\binom{n}{i} V(K[i],\mathbb{B}^n[n-i]),
\end{equation*}
where $\kappa_{m}$ denotes the $m$-dimensional volume of $\mathbb{B}^m$.

Associated with an $(n-1)$-tuple of bodies $K_2,\ldots, K_n \in \mathcal{K}^n$ is a finite Borel measure $S(K_2,\ldots,K_n,\cdot)$ on
$\mathbb{S}^{n-1}$, the {\it mixed area measure}, such that for all $K_1 \in \mathcal{K}^n\!$,
\begin{equation} \label{defmixedarea}
V(K_1,\ldots,K_n)=\frac{1}{n}\int_{\mathbb{S}^{n-1}} h(K_1,u)\,dS(K_2,\ldots,K_n,u).
\end{equation}

For $K \in \mathcal{K}^n$ and $0 \leq i \leq n - 1$, the measures $S_i(K,\cdot) := S(K[i],\mathbb{B}^n[n-1-i],\cdot)$ are called the \emph{area measures} of order $i$ of $K$. The measure $S_{n-1}(K,\cdot)$ is
also known as the \emph{surface area measure} of $K$. If $K$ has non-empty interior, then, by a theorem of Aleksandrov--Fenchel--Jessen (see, e.g.,
\textbf{\cite[\textnormal{p.\ 449}]{schneider93}}), each of the measures $S_i(K,\cdot)$, $1 \leq i \leq n - 1$, determines $K$ up to translations.
The centroid of every area measure of a convex body is at the origin, that is, for every $K \in \mathcal{K}^n$ and all $i \in \{0, \ldots, n - 1\}$,
\begin{equation*}
\int_{\mathbb{S}^{n-1}} u\,dS_i(K,u) = o.
\end{equation*}
\emph{Minkowski's existence theorem} states that a non-negative Borel measure $\mu$ on $\mathbb{S}^{n-1}$ is the
surface area measure of some $K \in \mathcal{K}^n$ with non-empty interior if and only if $\mu$ is not concentrated on a great subsphere of $\mathbb{S}^{n-1}$ and has centroid at the origin.

If $K \in \mathcal{K}^n$ has a $C^2$ support function, then each measure $S_i(K,\cdot)$, $0 \leq i \leq n - 1$, is absolutely continuous w.r.t.\ spherical Lebesgue measure. To make this more precise, let us recall the notion of mixed discriminants.
If $A_1, \ldots, A_m$ are symmetric real $k \times k$ matrices and $\lambda_1,\dots,\lambda_m \geq 0$, then
\begin{equation}\label{polydet}
\det(\lambda_1A_1+\cdots+\lambda_m A_m) = \sum_{j_1,\dots,j_k = 1}^m \mathrm{D}(A_{j_1},\ldots,A_{j_k})\lambda_{j_1}\cdots\lambda_{j_k},
\end{equation}
where the coefficients $\mathrm{D}(A_{j_1},\ldots,A_{j_k})$ are the {\it mixed discriminants} of $A_{j_1}, \ldots, A_{j_k}$ which depend only on $A_{j_1}, \ldots, A_{j_k}$ and are symmetric and multilinear in their arguments. Clearly, $\mathrm{D}(A,\dots,A) = \det(A)$ for any symmetric $k \times k$ matrix $A$. Moreover, $\mathrm{D}(B A_1,\ldots, B A_k) = \det(B)\mathrm{D}(A_1,\ldots,A_k)$,
\begin{equation} \label{basicMD}
\mathrm{D}(A,B,\ldots,B) = \frac{1}{n-1}\mathrm{tr}(\mathrm{cof}(B)A)
\end{equation}
for any symmetric $k \times k$ matrix $B$, and if $A_1,\ldots,A_k$ are positive semi-definite, then $\mathrm{D}(A_1,\dots,A_k) \geq 0$. Finally, if $K_1,\ldots,K_{n-1} \in \mathcal{K}^n$ have support functions $h_1,\ldots,h_{n-1} \in C^2(\mathbb{R}^n)$, then the density of $S(K_1,\ldots,K_{n-1},\cdot)$ is given by
\begin{equation} \label{densofmixedarea}
s(K_1,\ldots,K_{n-1},u) = \mathrm{D}(D^2h_1(u),\ldots,D^2h_{n-1}(u)), \quad u \in \mathbb{S}^{n-1}.
\end{equation}
In particular, for $K \in \mathcal{K}^n$ with support function $h \in C^2(\mathbb{R}^n)$, we have
\begin{equation} \label{densofsurfarea}
s_{n-1}(K,u) = \det D^2h(u),\quad u \in \mathbb{S}^{n-1}.
\end{equation}
Motivated by (\ref{densofmixedarea}) and (\ref{densofsurfarea}), we frequently use in subsequent sections the notation $s(h_1,\ldots,h_{n-1},\cdot)$, $s_{n-1}(h,\cdot)$, $\ldots$ instead of
$s(K_1,\ldots,K_{n-1},\cdot)$, $s_{n-1}(K,u)$, etc.

An origin-symmetric convex body $Z^{\mu} \in \mathcal{K}^n$ whose support function has an integral representation of the form
\begin{equation} \label{zonoidmeas}
h(Z^{\mu},u) = \int_{\mathbb{S}^{n-1}} |u \cdot v|\,d\mu(v), \qquad u \in \mathbb{S}^{n-1},
\end{equation}
with a (unique) even signed measure $\mu$ on $\mathbb{S}^{n-1}$ is called a \emph{generalized zonoid}. 
If $\mu$ is non-negative, then (\ref{zonoidmeas}) always defines a support function and the bodies obtained in this way are the origin-symmetric zonoids.
(see, e.g., \textbf{\cite[\textnormal{Chapter~3.5}]{schneider93}}).

\pagebreak

We turn now to the background material on spherical harmonics. To this end, let $\Delta_{\mathbb{S}}$ denote the spherical Laplacian on $\mathbb{S}^{n-1}$ and recall that it is
a second-order uniformly elliptic self-adjoint operator. We write $\mathcal{H}_k^n$ for the vector space of spherical harmonics of dimension $n$ and degree $k$ and denote its dimension by
\begin{equation} \label{nnk}
N(n,k) = \frac{n + 2k - 2}{n + k - 2} {n + k - 2 \choose n - 2} = \mathrm{O}(k^{n - 2}) \mbox{ as } k \rightarrow \infty.
\end{equation}
Spherical harmonics are (precisely) the eigenfunctions of $\Delta_{\mathbb{S}}$, more specific, for $Y_k \in \mathcal{H}_k^n$, we have
\begin{equation} \label{deltasmult}
\Delta_{\mathbb{S}} Y_k = -k(k + n - 2)\,Y_k.
\end{equation}

The spaces $\mathcal{H}_k^n$ are pairwise orthogonal subspaces of $L^2(\mathbb{S}^{n-1})$. Moreover, the Fourier series $f \sim \sum_{k=0}^{\infty} \pi_k f$
converges to $f$ in $L^2$ for every $f \in L^2(\mathbb{S}^{n-1})$, where $\pi_k: L^2(\mathbb{S}^{n-1}) \rightarrow \mathcal{H}_k^n$ denotes the orthogonal projection.
Letting $P_k^n \in C([-1,1])$ denote the \emph{Legendre polynomial} of dimension $n$ and degree $k$, we have
\begin{equation} \label{projleg}
(\pi_k f)(v) = \frac{N(n,k)}{\omega_n} \int_{\mathbb{S}^{n-1}}\!\! f(u)\, P_k^n(u\cdot v) \,du, \qquad v \in \mathbb{S}^{n-1},
\end{equation}
where $\omega_n$ denotes the surface area of $\mathbb{B}^n$ and integration is with respect to spherical Lebesgue measure. Since the orthogonal projection $\pi_k$ is self adjoint, it is consistent to extend it to the space $\mathcal{M}(\mathbb{S}^{n-1})$ of signed finite Borel measures by
\begin{equation*}
(\pi_k \mu)(v) = \frac{N(n,k)}{\omega_n}\int_{\mathbb{S}^{n-1}}\!\!P^n_k(u\cdot v)\, d\mu(u), \qquad v \in \mathbb{S}^{n-1}.
\end{equation*}
It can be shown easily that $\pi_k \mu \in \mathcal{H}_k^n$ for all $k \geq 0$ and that the formal Fourier series
$\mu \sim \sum_{k = 0}^\infty \pi_k f$ uniquely determines the measure $\mu$.

Throughout, we use $\bar{e} \in \mathbb{S}^{n-1}$ to denote a fixed but arbitrarily chosen pole of $\mathbb{S}^{n-1}$ and write $\mathrm{SO}(n-1)$ for the stabilizer
in $\mathrm{SO}(n)$ of $\bar{e}$. In the theory of spherical harmonics, a function or measure on $\mathbb{S}^{n-1}$ which is $\mathrm{SO}(n-1)$ invariant is often called \emph{zonal}.
Clearly, zonal functions depend only on the value of $u \cdot \bar{e}$. The subspace of zonal functions in $\mathcal{H}_k^n$ is $1$-dimensional for every $k \geq 0$ and spanned by $u \mapsto P_k^n(u \cdot \bar{e})$.
Since the spaces $\mathcal{H}_k^n$ are orthogonal, it is not difficult to show that any zonal measure $\mu \in \mathcal{M}(\mathbb{S}^{n-1})$
admits a series expansion of the form
\begin{equation} \label{expzonal}
\mu \sim \sum_{k=0}^{\infty} \frac{N(n,k)}{\omega_n}\, a_k^n[\mu]\,P_k^n(\,\,.\cdot \bar{e}),
\end{equation}
where
\begin{equation} \label{multleg}
a_k^n[\mu] = \omega_{n-1} \int_{-1}^1 P_k^n(t)\,(1-t^2)^{\frac{n-3}{2}}\,d\bar{\mu}(t).
\end{equation}
Here, we have used cylindrical coordinates $u = t\bar{e} + \sqrt{1 - t^2} v$ on $\mathbb{S}^{n-1}$ to identify the zonal measure $\mu$ with a measure $\bar{\mu}$ on $[-1,1]$. If $\mu$ is absolutely continuous with density $f$ w.r.t.\ spherical Lebesgue measure, we write $a_k^n[f\,]$ instead of $a_k^n[\mu]$.

\pagebreak

For the explicit computation of integrals of the form (\ref{multleg}) the following \emph{Formula of Rodrigues} for the Legendre polynomials
is often useful:
\begin{equation}\label{RodriguesF}
P_k^n(t)= (-1)^k \frac{\Gamma\left(\frac{n - 1}{2}\right)}{2^k \Gamma\left( \frac{n - 1}{2} + k\right)} (1 - t^2)^{-\frac{n-3}{2}}\left(\frac{d}{dt}\right)^k(1 - t^2)^{\frac{n-3}{2}+k}.
\end{equation}
It is clear from (\ref{RodriguesF}) that the derivative of a Legendre polynomial is itself a Legendre polynomial of higher dimension. Indeed,
\begin{equation}\label{derivLeg}
\frac{d}{dt} P_{k}^n(t) = \frac{k(k+n-2)}{n-1}P_{k-1}^{n+2}(t).
\end{equation}
There are several functional equations satisfied by the Legendre polynomials. One of the most noted is the following second order differential equation:
\begin{equation} \label{diffequpkn}
(1-t^2)\frac{d^{\,2}}{dt^2}P_k^n(t) -(n-1)t\frac{d}{dt}P_k^n(t) + k(k+n-2)P_k^n(t) = 0.
\end{equation}
In fact, this equation completely determines $P_k^n$ up to a constant factor.

Next, let us recall the well known \emph{Funk--Hecke Theorem}: If $\phi \in C([-1,1])$ and $\mathrm{T}_{\phi}: \mathcal{M}(\mathbb{S}^{n-1}) \rightarrow C(\mathbb{S}^{n-1})$ is defined by
\begin{equation} \label{Funkheckinttrafo}
(\mathrm{T}_{\phi}\mu)(u) = \int_{\mathbb{S}^{n-1}}\!\! \phi(u \cdot v)\,d\mu(v), \qquad u \in \mathbb{S}^{n-1},
\end{equation}
then the spherical harmonic expansion of $\mathrm{T}_{\phi}\mu$ is given by
\begin{equation} \label{funkhecke}
\mathrm{T}_{\phi}\mu \sim \sum_{k=0}^\infty a_k^n[\phi]\,\pi_k\mu,
\end{equation}
where the numbers $a_k^n[\phi]$ are given by (\ref{multleg}) and called the \emph{multipliers} of $\mathrm{T}_{\phi}$.

Integral transforms of the form (\ref{Funkheckinttrafo}) are closely related to the convolution between functions and measures on $\mathbb{S}^{n-1}$.
In order to recall its definition, first note that the convolution $\sigma \ast \tau$ of signed measures $\sigma, \tau$ on the compact Lie group $\mathrm{SO}(n)$ can be defined by
\[\int_{\mathrm{SO}(n)}\!\!\! f(\vartheta)\, d(\sigma \ast \tau)(\vartheta)=\int_{\mathrm{SO}(n)}\! \int_{\mathrm{SO}(n)}\!\!\! f(\eta \theta)\,d\sigma(\eta)\,d\tau(\theta), \qquad f \in C(\mathrm{SO}(n)).   \]
By identifying $\mathbb{S}^{n-1}$ with the homogeneous space $\mathrm{SO}(n)/\mathrm{SO}(n-1)$, one obtains a one-to-one correspondence of $C(\mathbb{S}^{n-1})$ and $\mathcal{M}(\mathbb{S}^{n-1})$, respectively, with right
$\mathrm{SO}(n-1)$ invariant functions and measures on $\mathrm{SO}(n)$, respectively. Using this correspondence, the convolution of measures on $\mathrm{SO}(n)$ induces a convolution product
on $\mathcal{M}(S^{n-1})$ (for more details see, e.g., \textbf{\cite{Schu06a}}). For this spherical convolution, zonal functions and measures play a particularly important role.
Let us therefore denote by $C(\mathbb{S}^{n-1},\bar{e})$ the set of continuous zonal functions on $\mathbb{S}^{n-1}$. Then, for $\mu \in \mathcal{M}(\mathbb{S}^{n-1})$, $f \in C(\mathbb{S}^{n-1},\bar{e})$, and $\eta \in \mathrm{SO}(n)$, we have
\begin{equation} \label{zonalconv}
(\mu \ast f)(\eta \bar{e}) =\int_{\mathbb{S}^{n-1}}\!\!f(\eta^{-1} u)\,d\mu(u).
\end{equation}

\pagebreak

Note that, by (\ref{zonalconv}), we have $(\vartheta \mu) \ast f = \vartheta(\mu \ast f)$ for every $\vartheta \in \mathrm{SO}(n)$,
where $\vartheta\mu$ is the image measure of $\mu$ under $\vartheta \in \mathrm{SO}(n)$. Moreover, from the identification of a
zonal function $f$ on $\mathbb{S}^{n-1}$ with a function $\bar{f}$ on $[-1,1]$, (\ref{multleg}), and (\ref{zonalconv}), we obtain
\begin{equation} \label{multconvpkn}
a_k^n[f] = \int_{\mathbb{S}^{n-1}}\!\!f(u) P_k^n(\bar{e}\cdot u)\,du
\end{equation}
and the Funk--Hecke Theorem implies that
\begin{equation} \label{convmulttransf}
\mu \ast f \sim \sum_{k=0}^{\infty} a_k^n[f]\,\pi_k\mu.
\end{equation}
Hence, convolution from the right induces a multiplier transformation. It is also easy to check from (\ref{zonalconv}) that the convolution of zonal functions and measures is Abelian and that for all $\mu, \tau \in \mathcal{M}(\mathbb{S}^{n-1})$ and every $f \in C(\mathbb{S}^{n-1},\bar{e})$,
\begin{equation} \label{convselfad}
\int_{\mathbb{S}^{n-1}}\!\! (\mu \ast f)(u)\,d\tau(u) = \int_{\mathbb{S}^{n-1}}\!\! (\tau \ast f)(u)\,d\mu(u).
\end{equation}

\vspace{0.1cm}

\begin{expls} \label{exps1} \end{expls}

\vspace{-0.2cm}

\begin{enumerate}
\item[(a)] The cosine transform $\mathrm{C}: \mathcal{M}(\mathbb{S}^{n-1}) \rightarrow C(\mathbb{S}^{n-1})$ is defined by
\begin{equation}
\mathrm{C}\mu(u) = \int_{\mathbb{S}^{n-1}}\!\!|u \cdot v|\,d\mu(v) = (\mu \ast |\bar{e}\cdot\,.\,|)(u), \qquad u \in \mathbb{S}^{n-1}.
\end{equation}
Using the Formula of Rodrigues, the multipliers $a_k^n[\mathrm{C}] := a_k^n[|\bar{e}\cdot\,.\,|]$ of the cosine transform can be easily computed to
\begin{equation} \label{multC}
a_k^n[\mathrm{C}] =  (-1)^\frac{n-2}{2}2\frac{1\cdot 3\cdots (k-3)}{(n+1)(n+3)\cdots (k+n-1)}
\end{equation}
for even $k$ and $a_k^n[\mathrm{C}]= 0$ for $k$ odd.

\item[(b)] Noting that $|\bar{e}\cdot\,.\,| = h([-\bar{e},\bar{e}],\cdot)$, we consider as a generalization of (a) integral transforms generated by generalized zonoids of revolution.
Indeed, for a (suitable) signed zonal measure $\mu$ on $\mathbb{S}^{n-1}$, let $\mathrm{T}_{\!Z^{\mu}}: \mathcal{M}(\mathbb{S}^{n-1}) \rightarrow C(\mathbb{S}^{n-1})$ be defined by
\begin{equation} \label{TZmu}
\mathrm{T}_{\!Z^{\mu}} \sigma = \sigma \ast h(Z^{\mu},\cdot) = \sigma \ast |\bar{e}\cdot\,.\,| \ast \mu = \mathrm{C}\sigma \ast \mu,
\end{equation}
where we have used (\ref{zonoidmeas}) and the commutativity of the convolution of zonal measures.
Hence, the multipliers $a_k^n[Z^{\mu}] := a_k^n[h(Z^{\mu},\cdot)]$ of $\mathrm{T}_{\!Z^{\mu}}$ are given by
\begin{equation}
a_k^n[Z^\mu] =  a_k^n[\mathrm{C}] a_k^n[\mu].
\end{equation}

\item[(c)] Generalizing now (b), we define for an arbitrary body of revolution $L \in \mathcal{K}^n$, the integral transform $\mathrm{T}_{\!L}: \mathcal{M}(\mathbb{S}^{n-1}) \rightarrow C(\mathbb{S}^{n-1})$ by
\[\mathrm{T}_L \sigma = \sigma \ast h(L,\cdot).\]
As in (b), we denote its multipliers by $a_k^n[L] := a_k^n[h(L,\cdot)]$.
\end{enumerate}

\pagebreak

In the final part of this section, we discuss the required material from functional analysis. First, let $H^s_e(\mathbb{S}^{n-1})$ denote the Sobolev space of \emph{even} functions on $\mathbb{S}^{n-1}$ with weak (covariant) derivatives up to order $s$ in $L^2(\mathbb{S}^{n-1})$ endowed with its standard norm $\|\cdot\|_{H^s}$. Then, by a classical result of Strichartz \textbf{\cite{Strichartz1981}},
\begin{equation} \label{strichartz}
\|f\|_{H^s}^2 = \sum_{k  = 0}^\infty (1+k^2)^s\|\pi_k f\|_{L^2}^2
\end{equation}
for $f \in H^s_e(\mathbb{S}^{n-1})$. Therefore, if a convolution transform $\mathrm{T}_{\mu}: L^2(\mathbb{S}^{n-1}) \rightarrow L^2(\mathbb{S}^{n-1})$ is generated by a zonal measure $\mu$, that is, $\mathrm{T}_{\mu}f = f \ast \mu$, satisfying
\begin{equation*}
a_k^n[\mu] = \mathrm{O}(k^{-\varrho}) \qquad \mbox{as } k \rightarrow \infty
\end{equation*}
for some integer $\varrho \geq 0$, then for a suitable constant $C_{n,\varrho}^s$,
\begin{equation}\label{smoothm}
\|\mathrm{T}_\mu f\|_{H^{s+\varrho}} \leq C_{n,\varrho}^s \|f\|_{H^{s}},
\end{equation}
for all $f \in H^s_e(\mathbb{S}^{n-1})$. Thus, if $f\in H_e^s(\mathbb{S}^{n-1})$, then $\mathrm{T}_\mu f \in H^{s+\varrho}_e(\mathbb{S}^{n-1})$.

We also recall a special case of the well known Sobolev embedding theorem for compact manifolds (see, e.g., \textbf{\cite[\textnormal{Chapter 2}]{Aubin1982}}). It states that the embedding
\begin{equation} \label{sobembed}
H^s_e(\mathbb{S}^{n-1}) \subseteq C^2(\mathbb{S}^{n-1})
\end{equation}
is continuous whenever $n < 2s - 3$.

Next, let $X$ and $Y$ be normed spaces and $U\subseteq X$ be an open subset of $X$. A function $f: U \to Y$ is called \emph{Fr\'echet differentiable} at $ x \in U$, if there exists a bounded linear operator $\mathrm{d}f(x) :X \to Y$, the \emph{Fr\'echet derivative} of $f$ at $x$, such that
\[\lim_{\|h\|_X \to 0} \frac{ \| f(x + h) - f(x) - \mathrm{d}f(x)h \|_{Y} }{ \|h\|_{X} } = 0.\]

Note that if $X$ and $Y$ are finite dimensional, then the Fr\'echet derivative of $f$ coincides with the usual derivative, represented in coordinates by the Jacobian of~$f$.
Moreover, Fr\'echet derivatives satisfy many of the basic properties of the usual derivative such as linearity and the product and chain rule.

\vspace{0.2cm}

\begin{expls} \label{exps2} \end{expls}
\vspace{-0.2cm}
\begin{enumerate}
\item[(a)] It was shown in \textbf{\cite{Martinez2001}} that the restriction of the cosine transform to continuous functions on $\mathbb{S}^{n-1}$
defines a bounded (and, hence, continuous) linear operator $\mathrm{C}: C(\mathbb{S}^{n-1}) \rightarrow C^2(\mathbb{S}^{n-1})$. This follows from explicit expressions
for the first and second (Fr\'echet) derivatives of $\mathrm{C}f$, when the latter is considered as a $1$-homogeneous function on $\mathbb{R}^n\backslash \{0\}$:
\begin{align*}
\mathrm{d}(\mathrm{C}f)(x) & = \int_{\mathbb{S}^{n-1}} \!\! \mathrm{sgn}(u \cdot x)uf(u)\,du, \\
\mathrm{d}^2(\mathrm{C}f)(x) & = \frac{2}{\|x\|} \int_{\mathbb{S}^{n-1} \cap\, x^{\bot}} \!\!\!\! u \otimes u f(u)\,du.
\end{align*}
\item[(b)] Suppose that $L \in \mathcal{K}^n$ is a body of revolution of class $C^2_+$. Restricting the integral transform $\mathrm{T}_L$, defined in Example \ref{exps1} (c),
to continuous functions on $\mathbb{S}^{n-1}$, defines a bounded linear operator $\mathrm{T}_L: C(\mathbb{S}^{n-1}) \rightarrow C^2(\mathbb{S}^{n-1})$. To see this, we compute the first and second derivatives of the $1$-homogeneous extension of $\mathrm{T}_Lf$ to $\mathbb{R}^n\backslash \{0\}$,
\begin{equation} \label{TLnotation}
\mathrm{T}_Lf(x) = \int_{\mathbb{S}^{n-1}}\!\! h(L(u),x)f(u)\,du.
\end{equation}
Here, $L(u)$ denotes the rotated copy of $L$ whose axis of symmetry is $u \in \mathbb{S}^{n-1}$. From (\ref{TLnotation}), it is immediate that
\begin{align*}
\mathrm{d}(\mathrm{T}_Lf)(x) & = \int_{\mathbb{S}^{n-1}}\!\! \nabla h(L(u),\cdot)(x) f(u)\,du, \\
\mathrm{d}^2(T_Lf)(x) & = \int_{\mathbb{S}^{n-1}}\!\!  D^2h(L(u),\cdot)(x) f(u)\,du.
\end{align*}

\item[(c)] Let $X_1, \ldots, X_n$, and $Y$ be normed spaces and let $\mathrm{T}: X_1\times \cdots \times X_n \to Y$ be a multilinear map which is bounded, that is,
\begin{equation*}
\|\mathrm{T}(x_1,\dots,x_n)\|_Y \leq C\|x_1\|_{X_1}\cdots \|x_n\|_{X_n}
\end{equation*}
for some $C > 0$ and all $x_1\in X_1,\ldots, x_n\in X_n$. Then $T$ is Fr\'echet differentiable on $X_1\times \cdots \times X_n $ and
\begin{equation*}
\mathrm{d}\mathrm{T}(x_1,\dots,x_n)(u_1,\dots, u_n) =  \sum_{k = 1}^{n}T(x_1,\dots,x_{k-1},u_k,x_{k+1},\dots,x_{n})
\end{equation*}
for all $x_1,u_1 \in X_1,\ldots, x_n, u_n\in X_n$.
\end{enumerate}

\vspace{0.2cm}

We require the following version of the inverse function theorem for Banach spaces, where we call a function $f$ defined on a subset $U$ of a normed space $X$ mapping to a normed space $Y$ \emph{continuously differentiable} if it is Fr\'echet differentiable on $U$ and $\mathrm{d}f$ is continuous as a function from $U\times X$ to $Y$.

\begin{theorem} \emph{(\!\!\textbf{\cite[\textnormal{Theorem 5.2.3}]{Hamilton1982}})} \label{invfctthm}
Let $X$ be a Banach space, $U \subseteq X$ an open subset, $Y$ a normed space, and suppose that $f: U \rightarrow Y$ is continuously differentiable. If $x\in U$ and $\mathrm{d}f(x)$ is invertible, then $f$ is a local diffeomorphism at $x$.
\end{theorem}

\vspace{0.8cm}

\centerline{\large{\bf{ \setcounter{abschnitt}{3}
\arabic{abschnitt}. Isoperimetric Inequalities}}}

\reseteqn \alpheqn \setcounter{theorem}{0}

\vspace{0.6cm}

In the following we first give the proof of Proposition 1. In the second part of this section we show that $V_{2}(\Phi_1 K)/V_{2}(K)$ is minimized
by Euclidean balls for \emph{every} continuous non-trivial Minkowski valuation $\Phi_1$ of degree $1$ compatible with rigid motions. This extends previous results from
\textbf{\cite{Schu06a}}, where an additional assumption of (weak) monotonicity was required.

\pagebreak

In order to prove Proposition 1, recall that for $1 \leq i \leq n - 1$ and $K, L \in \mathcal{K}^n$ having dimension at least $i + 1$, we have
\begin{equation} \label{aleksandrfench}
V(L,K[i],\mathbb{B}^n[n-i-1])^{i+1} \geq V(L[i+1],\mathbb{B}^n[n-i-1])V(K[i+1],\mathbb{B}^n[n-i-1])^{i}
\end{equation}
with equality if and only if $K$ and $L$ are homothetic. This is a well known consequence of the Aleksandrov--Fenchel inequality, where the equality cases are known
(see, e.g., \textbf{\cite[\textnormal{Chapter 7}]{schneider93}}).

\vspace{0.3cm}

\noindent {\it Proof of Proposition 1.} By a classical result of McMullen \textbf{\cite{McMullen77}}, the translation invariance and monotonicity of the Minkowski valuation $\Phi_i$ imply that it is continuous. Let $f \in L^1(\mathbb{S}^{n-1})$ be its (zonal) generating function.
Since we assume $\Phi_i$ to be non-trivial, there exists $K \in \mathcal{K}^n$ such that $\Phi_i K$ is not a singleton.
Consequently, the mean width $w(\Phi_iK)$ is strictly positive, that is,
\[0 < w(\Phi_i K) = \frac{2}{\omega_n}\int_{\mathbb{S}^{n-1}}\!\! h(\Phi_i K,u)\,du = \frac{2}{\omega_n}\int_{\mathbb{S}^{n-1}}\!\! (S_i(K,\cdot) \ast f)(u)\,du.   \]
Hence, by (\ref{convselfad}), (\ref{multconvpkn}), and the fact that $P_0^n(t) = 1$,
\[0 <\! \int_{\mathbb{S}^{n-1}}\!\! (S_i(K,\cdot) \ast f)(u)\,du =\! \int_{\mathbb{S}^{n-1}}\!\! (1 \ast f)(u)\,dS_i(K,u) = a_0^n[f]nV(K[i],\mathbb{B}^n[n-i]).  \]
Thus, $a_0^n[f] > 0$. Moreover, by homogeneity and translation invariance, $\Phi_i B(r,x)$ for an arbitrary ball $B(r,x)$ of radius $r > 0$ and center $x \in \mathbb{R}^n$, is given by
\begin{equation} \label{imageofball}
\Phi_i B(r,x) = r^i \Phi_i \mathbb{B}^n = r^i a_0^n[f] \mathbb{B}^n.
\end{equation}

In order to establish inequality (\ref{prop1inequ}), we first need to show that its right hand side is well defined, that is, $\mathrm{dim}\,\Phi_i K \geq i + 1$ for every $K \in \mathcal{K}^n$ whose dimension is at least $i + 1$. To this end, first assume that $i = n - 1$. Then, by (\ref{imageofball}) and our assumption of monotonicity of $\Phi_{n-1}$, it follows that $\Phi_{n-1}K$ has non-empty interior for every full dimensional $K \in \mathcal{K}^n$.

Let now $1 \leq i \leq n - 2$ and for a subspace $E \in \mathrm{Gr}_{i + 1,n}$, denote by $\mathbb{B}(E)$ the $i + 1$ dimensional unit ball in $E$. By the translation invariance, monotonicity, and the fact that $\Phi_i$ commutes with $\mathrm{SO}(n)$, it suffices to show that $\mathrm{dim}\,\Phi_i \mathbb{B}(E) \geq i + 1$ for \emph{some} $E \in \mathrm{Gr}_{i + 1,n}$.
By considering rotations that leave $E^{\bot}$ pointwise fixed, it follows that $\mathrm{p}_E\Phi_i \mathbb{B}(E) = r_{E,\Phi_i}\mathbb{B}(E)$ for some $r_{E,\Phi_i} \geq 0$, where $\mathrm{p}_E: \mathbb{R}^n \rightarrow E$ denotes the orthogonal projection. If $r_{E,\Phi_i} > 0$ for some $E \in \mathrm{Gr}_{i + 1,n}$, then the claim follows.

If, on the other hand, $\mathrm{p}_E\Phi_i \mathbb{B}(E) = \{o\}$ for all $E \in \mathrm{Gr}_{i + 1,n}$, then, by the monotonicity of $\Phi_i$, it follows that $\Phi_i K \subseteq E^{\bot}$ for every convex body $K \subseteq E$ and all $E \in \mathrm{Gr}_{i + 1,n}$. Moreover, by considering rotations that leave $E$ pointwise fixed, we have $\Phi_i K = r_{E,\Phi_i}^K\mathbb{B}(E^{\bot})$ for suitable
$r_{E,\Phi_i}^K \geq 0$ and every $K \subseteq E$ and $E \in \mathrm{Gr}_{i + 1,n}$. However, if $K \in \mathcal{K}^n$ contains the origin and $\mathrm{dim}\,K = i$, then, since $i \leq n - 2$, there exist \emph{distinct} $E, \bar{E} \in \mathrm{Gr}_{i + 1,n}$ such that $K \subseteq E$ and $K \subseteq \bar{E}$. Thus,
\[ \Phi_i K = r_{E,\Phi_i}^K\mathbb{B}(E^{\bot}) = r_{\bar{E},\Phi_i}^K\mathbb{B}(\bar{E}^{\bot}). \]
Clearly, this is possible only if $\Phi_i K = \{o\}$ for every $i$ dimensional $K \in \mathcal{K}^n$. But it was shown in \textbf{\cite[\textnormal{Lemma 2.5}]{schnschu}} that this implies that $\Phi_i$ is trivial, a contradiction. Hence, if $K$ has dimension at least $i + 1$, then so does $\Phi_iK$.

Next, note that by (\ref{convselfad}), we have
\begin{align*}
V(\Phi_i K,\Phi_iK[i],\mathbb{B}^n[n-i-1]) & =  \frac{1}{n} \int_{\mathbb{S}^{n-1}}\!\! h(\Phi_i K,u)\,dS_i(\Phi_iK,u) \\
& =  \frac{1}{n} \int_{\mathbb{S}^{n-1}}\!\! (S_i(\Phi_iK,\cdot) \ast f)(u)\,dS_i(K,u) \\
& = V(\Phi_i^2 K,K[i],\mathbb{B}^n[n-i-1])
\end{align*}
for every $K \in \mathcal{K}^n$. If $\mathrm{dim}\,K \geq i + 1$, then combining this with (\ref{aleksandrfench}), yields
\begin{align*}
\frac{V_{i+1}(\Phi_i K)^{i+1}}{V_{i+1}(K)^{i+1}} & = \frac{V(\Phi_i^2 K,K[i],\mathbb{B}^n[n-i-1])^{i+1}}{V(K[i+1],\mathbb{B}^n[n-i-1])^{i+1}} \\
&  \geq \frac{V(\Phi_i^2K[i+1],\mathbb{B}^n[n-i-1])}{{V(K[i+1],\mathbb{B}^n[n-i-1])}}=\frac{V_{i+1}(\Phi_i^2 K)}{V_{i+1}(K)}
\end{align*}
with equality if and only if $\Phi_i^2 K$ and $K$ are homothetic, which, by rearranging terms, proves (\ref{prop1inequ}).

Assume now that $f = h(L,\cdot)$ for some body of revolution $L \in \mathcal{K}^n$ of class $C^2_+$ and that balls are the only solutions to the fixed-point problem $\Phi_i^2 K = \alpha K$.
Then, by the first part of the proof, only for balls $\psi_i(K) = V_{i+1}(\Phi_iK)/V_{i+1}(K)^i$ can attain a minimum among bodies $K \in \mathcal{K}^n$ of dimension at least $i + 1$. Hence, it remains to show that $\psi_i$ actually attains a minimum on this set. To this end, let $m_L = \min_{v \in \mathbb{S}^{n-1}}h(L,v)$ and note that $m_L > 0$, since $L$ is of class $C^2_+$, and $m_L = m_{L(u)}$ for every $u \in \mathbb{S}^{n-1}$. Hence, by (\ref{TLnotation}),
\[h(\Phi_{i}K,u) = \mathrm{T}_L S_{i}(K,\cdot)(u) = \int_{\mathbb{S}^{n-1}}\!\! h(L(u),v)\,dS_i(K,v) \geq m_L nV(K[i],\mathbb{B}^{n}[n-i]) \]
for every $u \in \mathbb{S}^{n-1}$ and all $K \in \mathcal{K}^n$ of dimension at least $i+1$. Consequently, we have $\Phi_iK \supseteq m_L nV(K[i],\mathbb{B}^{n}[n-i]) \mathbb{B}^n$ and, thus,
there exists $c_{n,i} > 0$ such that
\begin{equation} \label{proofprop117}
\psi_i(K) \geq  m_L^{i+1}c_{n,i} \frac{V_i(K)^{i+1}}{V_{i+1}(K)^{i}}
\end{equation}
for all $K \in \mathcal{K}^n$ of dimension at least $i+1$. For $2 \leq i \leq n - 1$, it was proved in \textbf{\cite{Gritzmann1988}}, and for $i = 1$, it is a classical fact (cf.\ also \textbf{\cite{Gritzmann1988}}) that there exists $\bar{c}_{n,i} > 0$ such that
\begin{equation}\label{ineq1}
V_i(K)^i\geq \bar{c}_{n,i} D(K)V_{i+1}(K)^{i-1}
\end{equation}
for all $K \in \mathcal{K}^n$, where $D(K)$ denotes the diameter of $K$. Combining (\ref{proofprop117}) and (\ref{ineq1}) with the isoperimetric inequality between consecutive intrinsic volumes (that is, the special case $L = \mathbb{B}^n$ of (\ref{aleksandrfench})), now yields $b_{n,i} > 0$ such that
\begin{equation} \label{ineq3}
\psi_i(K)\geq  m_L^{i+1}b_{n,i} D(K) V_{i+1}(K)^{-\frac{1}{i+1}}
\end{equation}
for all $K \in \mathcal{K}^n$ with $\mathrm{dim}\,K \geq i + 1$.

\pagebreak

Let $K_j \in \mathcal{K}^n$ be a sequence of bodies of dimension at least $i + 1$ such that $\psi_i(K_j) \to \inf \{\psi_i(K): K\in  \mathcal{K}^n, \mathrm{dim}\,K \geq i + 1\}$. By the translation and scaling invariance of $\psi_i$, we may assume that $K_j$ contains the origin and that $D(K_j) = 1$ for every $j \in \mathbb{N}$. Hence, by Blaschke's selection theorem (see, e.g., \textbf{\cite[\textnormal{Section 1.8}]{schneider93}}) $K_j$ admits a convergent subsequence. For simplicity let us assume this sequence to be $K_j$ itself and let $\bar{K}$ denote its limit.
If $\mathrm{dim}\,\bar{K} < i + 1$, then $V_{i+1}(\bar{K}) = 0$ and, hence, $\lim_{j\to \infty} V_{i+1}(K_j) = 0$. Thus, by (\ref{ineq3}), we have
\[\inf_{\mathrm{dim}\,K \geq i + 1}\psi_i(K) = \lim_{j\to \infty} \psi_i(K_j)\geq m_L^{i+1}b_{n,i} \lim_{j\to \infty} V_{i+1}(K_j)^{-\frac{1}{i+1}} = \infty,\]
which is clearly a contradiction. Hence, $\bar{K}$ has dimension at least $i+1$ and, thus, $\psi_i$ attains a minimum at $\bar{K}$. \hfill $\blacksquare$

\vspace{0.3cm}

Let us make two remarks on the assumptions of Proposition 1. First note that if the monotone Minkowski valuation $\Phi_i$ is generated by a function $f \in L^1(\mathbb{S}^{n-1})$ such that $f \geq m_f > 0$, then the statement still holds true, just replace $m_L$ in the proof by $m_f$. (This part of the proof was a refinement of arguments taken from \textbf{\cite{Hofstaetteretal2022}}.) We have chosen to state Proposition 1 in the introduction with the assumption of $C^2_+$ regularity as this is the novel aspect of our proposed new approach towards Petty's conjecture. Second, note that the assumption of monotonicity was used only to conclude that $\Phi_i$ maps bodies of dimension at least $i + 1$ to such bodies. It is an open problem whether this is true for all merely continuous non-trivial Minkowski valuations intertwining rigid motions. However, since a Minkowski valuation $\Phi_i$ generated by a (non-zero) convex body $L \in \mathcal{K}^n$, satisfies
\[h(\Phi_{i}K,v) = \mathrm{T}_L S_{i}(K,\cdot)(v) = nV(L(v),K[i],\mathbb{B}^{n}[n-i-1])\]
for every $K \in \mathcal{K}^n$ and $v \in \mathbb{S}^{n-1}$, it follows that such $\Phi_i$ is not only monotone but maps bodies of dimension $i + 1$ to bodies with non-empty interior.

While the assumption of monotonicity was critical in our proof of Proposition 1 and in the determination of minimizers of $V_{2}(\Phi_1 K)/V_{2}(K)$ in \textbf{\cite{Schu06a}}, we will show in the second part of this section, that for the latter it can be dropped altogether. To this end, recall that the second order differential operator $\Box_n$, defined by
\[\Box_n h = h +  \frac{1}{n-1}\Delta_{\mathbb{S}}h \]
for $h \in C^2(\mathbb{S}^{n-1})$, relates the support function $h(K,\cdot)$ of a convex body $K \in \mathcal{K}^n$ with its first-order area measure $S_1(K,\cdot)$ by
\begin{equation} \label{boxhks1}
\Box_n h(K,\cdot) = S_1(K,\cdot),
\end{equation}
where (\ref{boxhks1}) has to be understood in a distributional sense if $h(K,\cdot)$ is not $C^2$. From the definition of $\Box_n$ and (\ref{deltasmult}), we see
that for $h \in C(\mathbb{S}^{n-1})$ and every $k \geq 0$,
\begin{equation} \label{boxnmult}
\pi_k \Box_n h = \frac{(1-k)(k+n-1)}{n-1} \pi_k h.
\end{equation}
In particular, $\Box_n$ acts as a multiplier transformation and since such operators clearly commute, we note for later quick reference that, by (\ref{convmulttransf}) and (\ref{boxnmult}), we have
\begin{equation} \label{boxconvcomm}
\mathrm{T}_{f}\,\Box_n = \Box_n \mathrm{T}_f.
\end{equation}
After these preparations, we are now in a position to determine the minimizers of $V_{2}(\Phi_1 K)/V_{2}(K)$ without any monotonicity assumption.

\begin{theorem} \label{mainthmdeg1} Let $\Phi_1: \mathcal{K}^n \rightarrow \mathcal{K}^n$ be a non-trivial continuous translation invariant Minkowski valuation of degree $1$ which commutes with $\mathrm{SO}(n)$.
If $K \in \mathcal{K}^n$, then
\begin{equation} \label{vader17}
V_2(\Phi_1K) \geq \frac{V_2(\Phi_1\mathbb{B}^n)}{V_2(\mathbb{B}^n)}\,V_2(K).
\end{equation}
\end{theorem}

\noindent {\it Proof.} Let $f \in \mathbb{S}^{n-1}$ be the generating function of $\Phi_1$. As in the first part of the proof of Proposition 1, it follows that $a_0^n[f] > 0$, since $\Phi_1$ is non-trivial. Thus, by (\ref{imageofball}), we have $V_2(\Phi_1\mathbb{B}^n)/V_2(\mathbb{B}^n) =  a_0^n[f]^2$. Consequently, we want to show that $V_2(\Phi_1K) - a_0^n[f]^2 V_2(K) \geq 0$ for all $K \in \mathcal{K}^n$ or, equivalently, that
\begin{equation} \label{toshow17}
V(\Phi_1K[2],\mathbb{B}^n[n-2]) - a_0^n[f]^{2\,} V(K[2],\mathbb{B}^n[n-2]) \geq 0.
\end{equation}
To see this, note that, by (\ref{defmixedarea}), (\ref{toshow17}) is equivalent to
\[\int_{\mathbb{S}^{n-1}}\!\! h(\Phi_1K,u)\,dS_1(\Phi_1K,u) - a_0^n[f]^2\int_{\mathbb{S}^{n-1}}\!\! h(K,u)\,dS_1(K,u) \geq 0 \]
which, by (\ref{boxhks1}) and Theorem \ref{thmdorrek}, in turn is equivalent to
\[\int_{\mathbb{S}^{n-1}}\!\! \mathrm{T}_f\Box_nh(K,\cdot)\,\Box_n\mathrm{T}_f\Box_nh(K,\cdot) - a_0^n[f]^2\int_{\mathbb{S}^{n-1}}\!\! h(K,\cdot)\,\Box_nh(K,\cdot)  \geq 0. \]
Using now (\ref{boxconvcomm}), it follows that we need to show that
\[\int_{\mathbb{S}^{n-1}}\!\! h(K,\cdot) \Box_n((\mathrm{T}_f\Box_n)^2  - a_0^n[f]^2\mathrm{Id})\,h(K,\cdot) \geq 0. \]
holds for all $K \in \mathcal{K}^n$. Clearly, this follows if we prove that the multiplier transform $\mathrm{E}=\Box_n((\mathrm{T}_f\Box_n)^2  - a_0^n[f]^2\mathrm{Id})$ is positive semi-definite, which means that
all its multipliers are non-negative. To see this, note that, by (\ref{boxnmult}), $a_0^n[\mathrm{E}] = 0$, $a_1^n[\mathrm{E}] = 0$  and, $a_k^n[\Box_n] < 0$ for all $k \geq 2$. Hence, it remains to show that
\begin{equation} \label{Grogu}
|a_k^n[\mathrm{T}_f\Box_n]| \leq a_0^n[f]
\end{equation}
holds for all $k \geq 2$. But this was proved by Schneider in \textbf{\cite[\textnormal{Lemma 4.12}]{Schneider1974c}}. \hfill $\blacksquare$

\vspace{0.3cm}

The proof of Theorem \ref{mainthmdeg1} shows that equality in (\ref{vader17}) among bodies $K \in \mathcal{K}^n$ with $\mathrm{dim}\,K \geq 2$ holds precisely for balls as long as inequality (\ref{Grogu}) is strict
for all $k \geq 2$. Kiderlen \textbf{\cite{kiderlen05}} showed that this is the case if $\Phi_1$ is monotone.

\pagebreak

The above proof also allows to show the following stronger version of inequality (\ref{vader17}) if $\Phi_1$ is generated by a convex body $L \in \mathcal{K}^n$ (first established by Saroglou and Zvavitch \textbf{\cite{SaroglouZvavitch}} for $\Pi_1$),
\begin{equation*}
V(\Phi_{1\!}K[2],\mathbb{B}^n[n-2])_{\!} \geq  \frac{a_0^n[f]^2}{(n-1)^2}V(K[2],\mathbb{B}^n[n-2]) + \frac{n(n-2)a_0^n[f]^2}{(n-1)^2\omega_n}V(K,\mathbb{B}^n[n-1])^2
\end{equation*}
for all $K\in \mathcal{K}^n$. This is an immediate consequence of Theorem \ref{spectralgap} below which implies that the operator
\[\mathrm{\bar{E}} = -\Box_n((n-1)^2(\mathrm{T}_L\Box_n)^2 - a_0^n[f]^2\mathrm{Id} - n(n-2)a_0^n[f]^2 \pi_0) \]
is positive semi-definite.

\vspace{1cm}

\centerline{\large{\bf{ \setcounter{abschnitt}{4}
\arabic{abschnitt}. Auxiliary Results}}}

\reseteqn \alpheqn \setcounter{theorem}{0}

\vspace{0.6cm}

In the following we denote by $\mathcal{S}^{n,2}_+$ the subset of $C^2(\mathbb{S}^{n-1})$ consisting of all $h$ for which $D^2h>0$ or, equivalently, $\mathcal{S}^{n,2}_+$ is the set of all support functions of convex bodies in $\mathcal{K}^n$ of class $C^2_+$. It is not difficult to see that $\mathcal{S}^{n,2}_+$ is an open convex cone in $C^2(\mathbb{S}^{n-1})$. Our goal in this section is to compute the Fr\'echet derivatives of iterations $\Phi^m_i$ (that is, the composition of $\Phi_i$ with itself $m$ times) for sufficiently regular Minkowski valuations $\Phi_i$ of degree $i$ intertwining rigid motions when they are considered as operators on $\mathcal{S}^{n,2}_+$,
\begin{equation*}
        \Phi_i(h) = s_i(h,\cdot) \ast f, \qquad h \in \mathcal{S}^{n,2}_+,
\end{equation*}
where $f \in L^1(\mathbb{S}^{n-1})$ denotes the generating function of $\Phi_i$. This will turn out to be easy as soon as we compute the Fr\'echet derivatives of the area measure densities $s_i(h,\cdot)$ on $\mathcal{S}^{n,2}_+$. The Fr\'echet derivatives of $\Phi^m_i$ then follow from basic properties of derivatives: linearity and the chain rule. A first step towards this is contained in the following proposition.

\begin{prop}\label{Sdiff} The (multilinear extension of the) mixed area measure density map $s: C^2(\mathbb{S}^{n-1}) \times \cdots \times C^2(\mathbb{S}^{n-1}) \rightarrow C(\mathbb{S}^{n-1})$, given by
\[s(h_1,\ldots,h_{n-1},\cdot) = \mathrm{D}(D^2h_1,\ldots,D^2h_{n-1}),    \]
is Fr\'echet differentiable with derivative
\begin{equation*}
\mathrm{d}s(h_1,\dots,h_{n-1},\cdot)(g_1,\dots,g_{n-1})  =  \sum_{k = 1}^{n-1}s(h_1,\dots,h_{k-1},g_k,h_{k+1},\dots,h_{n-1},\cdot)
\end{equation*}
for $h_1,\ldots,h_{n-1}, g_1,\ldots, g_{n-1} \in C^2(\mathbb{S}^{n-1})$.
\end{prop}

\pagebreak

\noindent {\it Proof.} Since the mixed discriminant $\mathrm{D}$ is a bounded multilinear map with respect to a fixed norm (say the Frobenius norm $\Vert \cdot \Vert_{\mathrm{F}})$, there exists a $c > 0$ such that
\begin{equation}\label{bm}
|\mathrm{D}(A_1,\dots,A_{n-1})| \leq  c\, \Vert A_1 \Vert_{\mathrm{F}} \cdots \Vert A_{n-1} \Vert_\mathrm{F}
\end{equation}
for all $(n - 1) \times (n - 1)$  matrices $A_1,\dots, A_{n-1}$. Substituting $D^2h_1(u),\dots,D^2h_{n-1}(u)$ in (\ref{bm}) for $u \in \mathbb{S}^{n-1}$, we obtain
\begin{equation}\label{pointwise}
|s(h_1,h_2,\dots,h_{n-1},u)| \leq  c\, \Vert D^2h_1(u) \Vert_{\mathrm{F}} \cdots \Vert D^2h_{n-1}(u) \Vert_{\mathrm{F}}.
\end{equation}
Taking the supremum in (\ref{pointwise}) and noting that there exists $\bar{c} > 0$ such that
\begin{equation*}
\sup_{u\in\mathbb{S}^{n-1}}\Vert D^2h(u)\Vert_{\mathrm{F}} \leq \bar{c}\, \Vert h \Vert_{C^2(\mathbb{S}^{n-1})}
\end{equation*}
for all $h\in C^2(\mathbb{S}^{n-1})$, we conclude that $s$ is a bounded multilinear map from $C^2(\mathbb{S}^{n-1}) \times \cdots \times C^2(\mathbb{S}^{n-1})$ to $C(\mathbb{S}^{n-1})$. The claim now follows from an application of Example \ref{exps2} (c).
\hfill $\blacksquare$

\vspace{0.3cm}

For the following immediate consequence of Proposition \ref{Sdiff}, we write $h_{\mathbb{B}^n}$ for the support function of $\mathbb{B}^n$.

\begin{koro}\label{diffsurf} For $1 \leq i \leq n - 1$, the (extension of the) area measure density map $s_i: C^2(\mathbb{S}^{n-1}) \rightarrow C(\mathbb{S}^{n-1})$, given by
\[s_i(h,\cdot) = s(h[i],h_{\mathbb{B}^n}[n-i-1],\cdot) = \mathrm{D}(D^2h[i],\mathrm{Id}[n-1-i]),    \]
is Fr\'echet differentiable with derivative
\[\mathrm{d}s_i(h,\cdot)g = i \mathrm{D}(D^2g,D^2h[i-1],\mathrm{Id}[n-1-i])  \]
for all $h, g\in C^2(\mathbb{S}^{n-1})$.
\end{koro}
\noindent {\it Proof.} Define $\mathrm{J}: C^2(\mathbb{S}^{n-1}) \to C^2(\mathbb{S}^{n-1}) \times \cdots \times C^2(\mathbb{S}^{n-1})$, by $\mathrm{J}h = (h[i],h_{\mathbb{B}^n}[n-i-1])$.
Then, clearly, $\mathrm{J}$ is Fr\'echet differentiable with derivative $\mathrm{d}\mathrm{J}(h)g = (g[i],0[n-i-1])$ for all $h, g\in C^2(\mathbb{S}^{n-1})$.
Since $s_i = s \circ \mathrm{J}$, the chain rule implies that
\[\mathrm{d}s_i(h,\cdot)g = \mathrm{d}s(\mathrm{J}h,\cdot)\mathrm{d}\mathrm{J}(h)g = \mathrm{d}s(h[i],h_{\mathbb{B}^n}[n-i-1],\cdot)(g[i],0[n-i-1]).\]
Thus, by Proposition \ref{Sdiff},
\[\mathrm{d}s_i(h,\cdot)g  =  is(g,h[i-1],h_{\mathbb{B}^n}[n-i-1],\cdot) = i \mathrm{D}(D^2g,D^2h[i-1],\mathrm{Id}[n-1-i])\]
for all $h, g\in C^2(\mathrm{S}^{n-1})$. \hfill $\blacksquare$

\vspace{0.3cm}

We are now in a position to compute the Fr\'echet derivative of $\Phi_i^m$ on $\mathcal{S}^{n,2}_+$ for every $m \geq 1$. In order to simplify expressions, we sometimes normalize $\Phi_i$ such that $\Phi_i \mathbb{B}^n = \mathbb{B}^n$ or, equivalently, $a_0^n[f] = 1$ for its generating function $f \in L^1(\mathbb{S}^{n-1})$.

\begin{koro}\label{diffval} Let $1 \leq i \leq n - 1$ and $\Phi_i: \mathcal{K}^n \rightarrow \mathcal{K}^n$ be a continuous translation invariant Minkowski of degree $i$ which commutes with $\mathrm{SO}(n)$. If the convolution transform $\mathrm{T}_f: C(\mathbb{S}^{n-1}) \rightarrow C^2(\mathbb{S}^{n-1})$ is bounded for the generating function $f \in L^1(\mathbb{S}^{n-1})$ of $\Phi_i$, then $\Phi_i^m$ is Fr\'echet differentiable on $\mathcal{S}^{n,2}_+$ for every $m \geq 1$ with derivative given for $m = 1$, by
\begin{equation}\label{pval1}
\mathrm{d}\Phi_i(h)g = i\,\mathrm{T}_f \mathrm{D}(D^2g,D^2h[i-1],\mathrm{Id}[n-1-i]),
\end{equation}
and for $m > 1$, recursively by
\begin{equation}  \label{crval}
\mathrm{d}\Phi_i^m(h)g = \mathrm{d}\Phi_i(\Phi_i(h))\mathrm{d}\Phi_i^{m-1}(h)g
\end{equation}
for all $h\in \mathcal{S}^{n,2}_+$ and $g\in C^2(\mathbb{S}^{n-1})$. In particular, if $\Phi_i \mathbb{B}^n = \mathbb{B}^n$, then
\begin{equation}\label{diffphicompeqb}
\mathrm{d}\Phi_i^m(h_{\mathbb{B}^n})g = (i\,\Box_n \mathrm{T}_f)^m g
\end{equation}
for all $m \geq 1$ and $g\in C^2(\mathbb{S}^{n-1})$.
\end{koro}
\noindent {\it Proof.} In order to see (\ref{pval1}), note that $\Phi_i = \mathrm{T}_f \circ s_i$. Since we assume $\mathrm{T}_f$ to be a bounded linear map, we have
$\mathrm{d}\mathrm{T}_f(h)g = \mathrm{T}_{f\,}g$ for all $h \in C(\mathbb{S}^{n-1}), g \in C^2(\mathbb{S}^{n-1})$. Therefore, (\ref{pval1}) follows from the chain rule and Corollary \ref{diffsurf}.

The recursive expression (\ref{crval}) for $\mathrm{d}\Phi_i^m$ is a direct consequence of the chain rule. Finally, to see (\ref{diffphicompeqb}), we evaluate (\ref{pval1}) at $h_{\mathbb{B}^n}$ and use (\ref{basicMD}) and the fact that $\mathrm{tr}\,D^2g = (n - 1)\Box_ng$, to obtain
\begin{equation*}
\mathrm{d}\Phi_i(h_{\mathbb{B}^n})g  = i  \mathrm{T}_{f\,}\Box_n g = i\,\Box_n \mathrm{T}_{f\,} g
\end{equation*}
for all $g\in C^2(\mathbb{S}^{n-1})$ (the second equality is just (\ref{boxconvcomm})). Hence, by our assumption that $\Phi_i h_{\mathbb{B}^n} = h_{\mathbb{B}^n}$, (\ref{crval}) yields
\[\mathrm{d}\Phi_i^m(h_{\mathbb{B}^n})g = \mathrm{d}\Phi_i(h_{\mathbb{B}^n})\mathrm{d}\Phi_i^{m-1}(h_{\mathbb{B}^n})g = i\,\Box_n \mathrm{T}_{f\,} \mathrm{d}\Phi_i^{m-1}(h_{\mathbb{B}^n})g,\]
for all $g\in C^2(\mathbb{S}^{n-1})$ and (\ref{diffphicompeqb}) follows by induction. \hfill $\blacksquare$

\vspace{1cm}

\centerline{\large{\bf{ \setcounter{abschnitt}{5}
\arabic{abschnitt}. Spectral Gap}}}

\reseteqn \alpheqn \setcounter{theorem}{0}

\vspace{0.6cm}

In this section we conclude our preparations for the proof of our main result, by establishing a new spectral gap for the multiplier transform $\mathrm{T}_L$ generated by an origin-symmetric body of revolution $L \in \mathcal{K}^n$. To this end, we first collect some classical facts about the relative extremals of Legendre polynomials.

\begin{prop} \label{succrelmaxpkn} \textnormal{(\hspace{-0.3cm}\,\textbf{\cite[\textnormal{Section 7.8}]{Szegoe1975}})} Let $\nu_k^n[1], \ldots, \nu_k^n[[\frac{k}{2}]]$ denote the successive relative maxima of $|P_k^n(t)|$ as $t$ decreases from $1$ to $0$. Then the following holds:
\begin{itemize}
\item[(i)] $1 > \nu_k^n[1] > \nu_k^n[2] > \cdots > \nu_k^n[[\frac{k}{2}]]$ for all $k \geq 2$;
\item[(ii)] for every $r \geq 1$, we have $\nu_k^n[r] > \nu_{k+1}^n[r]$ for all $k \geq r + 1$.
\end{itemize}
\end{prop}

\pagebreak

As a first simple consequence of Proposition \ref{succrelmaxpkn}, we note the following (probably well known) universal lower bound for Legendre polynomials.

\begin{koro} \label{absmin} If $k \geq 2$ is even, then for every $t \in [-1,1]$,
\[P_k^n(t) \geq -\frac{1}{n-1}. \]
\end{koro}
\noindent {\it Proof.} By the Formula of Rodrigues (\ref{RodriguesF}), $P_2^n$ is given by
\begin{equation} \label{p2nexpl}
P_2^n(t) = \frac{nt^2 -1}{n-1}.
\end{equation}
Thus, an elementary computation yields
$$
\nu^n_2[1] = \frac{1}{n-1}.
$$
Hence, by Proposition \ref{succrelmaxpkn} we have for $k \geq 2$,
\begin{equation} \label{decrmin}
\frac{1}{n-1} = \nu^n_2[1] > \nu^n_k[1] > \nu^n_k[r],
\end{equation}
for all $r = 2,\dots,\left[\frac{k}{2}\right]$. In particular, if $k$ is even, then (\ref{decrmin}) and the fact that the relative extremals of $P_k^n$ alternate signs imply that  $-\nu^n_k[1]$ is the global minimum of $P_k^n$ on the real line. This together with (\ref{decrmin}) yields the desired inequality.  \hfill $\blacksquare$

\vspace{0.3cm}

The next lemma provides explicit conditions for a sufficiently regular function on $[-1,1]$ to give rise to a support function of a convex body of revolution in $\mathbb{R}^n$.

\begin{lem} \label{Cphi} Suppose that $\phi \in C^2([-1,1])$. Then $h \in C(\mathbb{R}^n)$, defined by
\[h(x) = \left \{ \begin{array}{ll} \|x\|\, \phi \hspace{-0.1cm} \left(\! \displaystyle{\frac{x \cdot \bar{e}}{\|x\|} } \right ) & \mbox{ for all } x \neq o, \\ 0 & \mbox{ for } x = o,     \end{array}    \right .  \]
is the support function of a convex body of revolution $K_{\phi} \in \mathcal{K}^n$ if and only if
\[\phi(t) - t\phi'(t) \geq 0 \qquad \mbox{and} \qquad (1-t^2)\phi''(t) + \phi(t) - t\phi'(t) \geq 0  \]
for all $t \in [-1,1]$. Moreover, $K_{\phi}$ is of class $C^2_+$ if and only if these inequalities are strict for all $t \in [-1,1]$.
\end{lem}

\noindent {\it Proof.} We first show that for $u \in \mathbb{S}^{n-1}$,
\begin{equation} \label{hessian1717}
D^2h (u) =  \left(\phi(u \cdot \bar{e})-(u\cdot \bar{e})\phi'(u\cdot \bar{e})\right)\mathrm{p}_{u^\bot} + \phi''(u\cdot \bar{e}) (\mathrm{p}_{u^\bot}\bar{e} \otimes \mathrm{p}_{u^\bot}\bar{e}),
\end{equation}
where $\mathrm{p}_{u^\bot}= \mathrm{Id} - u\otimes u$ denotes the orthogonal projection onto the hyperplane $u^\bot$. To see this, we compute the first partial derivatives of $h$ at a non-zero $x \in \mathbb{R}^n$,
\begin{equation}\label{partial1}
\frac{\partial h}{\partial x_i}(x)  = \phi\left(\frac{x\cdot \bar{e}}{\|x\|}\right) \frac{x_i}{\|x\|} + \phi'\left(\frac{x\cdot \bar{e}}{\|x\|}\right)\left(e_i - \frac{(x\cdot \bar{e})x_i}{\|x\|^2}\right).
\end{equation}

\pagebreak

Differentiating (\ref{partial1}) to find the second order derivatives, yields
\begin{multline*}
\frac{\partial^2 h}{\partial x_i\partial x_j}(x)  =   \left[\phi\left(\frac{x\cdot \bar{e}}{\|x\|}\right)- \frac{x\cdot \bar{e}}{\|x\|}\phi'\left(\frac{x\cdot \bar{e}}{\|x\|}\right)\right]\left(\frac{\delta_{ij}}{\|x\|}-\frac{x_i x_j}{\|x\|^3}\right)\\
+ \phi''\left(\frac{x\cdot \bar{e}}{\|x\|}\right)\left( \frac{e_j}{\|x\|} - \frac{(x\cdot \bar{e})x_j}{\|x\|^3}\right)\left( e_i - \frac{(x\cdot \bar{e})x_i}{\|x\|^2}\right)
\end{multline*}
and, hence, we obtain the desired result.

By (\ref{hessian1717}), the hessian $D^2h$ at $u \in \mathbb{S}^{n-1}$ can be written as
$$
D^2h(u)= g_1(u\cdot \bar{e}) (\mathrm{Id} - u\otimes u - v\otimes v) + g_2(u\cdot \bar{e})v\otimes v,
$$
where $g_1(t)  = \phi(t)-t\phi'(t)$, $g_2(t)  = (1-t^2)\phi''(t) + \phi(t) -t\phi'(t)$ for all $t\in [-1,1]$, and
$$v = \frac{\mathrm{p}_{u^\bot}\bar{e}}{\|\mathrm{p}_{u^\top}\bar{e}\|} = \frac{\mathrm{p}_{u^\bot}\bar{e}}{\sqrt{1-(u\cdot \bar{e})^2}}.$$
This yields an explicit spectral decomposition of $D^2h(u)$, where the eigenvalues are $g_1(u\cdot \bar{e})$ with multiplicity $n-2$,  $g_2(u\cdot \bar{e})$  and $0$  with multiplicity $1$. The eigenspace of the eigenvalue $0$ is the line spanned by $u$ corresponding to the fact that $D^2h(u)$ is the Hessian of a homogeneous function and, thus, orthogonal to $u$. \hfill $\blacksquare$

\vspace{0.3cm}

The following consequence of Lemma \ref{Cphi} is crucial for the proof of our spectral gap result.

\begin{prop} \label{intervals1742} Let $k \geq 2$ and $I_k^n$, $J_k^n \subseteq \mathbb{R}$ denote the intervals of all $\lambda$ and $\gamma$, respectively, for which
\[h_{\lambda}(u) = 1 + \lambda P_k^n(\bar{e} \cdot u) \qquad \mbox{and} \qquad s_{\gamma}(u) = 1 + \gamma P_k^n(\bar{e} \cdot u)    \]
are the support function of a convex body $K_{\lambda} \in \mathcal{K}^n$ and the density of the surface area measure of a convex body $K_{\gamma} \in \mathcal{K}^n$, respectively. Then
\[I_k^n \subseteq \left [-\frac{1}{(k(n+k-2)-1)\nu_k^n[1]},\frac{n-1}{(k-1)(n+k-1)}  \right ] \quad \mbox{ and } \quad J_k^n = \left [ -1,\frac{1}{\nu_k^n[1]}  \right ]  \]
and the inclusion for $I_k^n$ becomes equality for $k = 2$ and $\nu_2^n[1] = \frac{1}{n-1}$. Moreover, for $k = 2$, $K_\lambda$ is of class $C_2^+$ if only if $\lambda \in \mathrm{int}\,I_2^n$.
\end{prop}

\noindent {\it Proof.} In order to simply notation, we write in the following $P$ instead of $P_k^n$. According to Lemma \ref{Cphi}, $\lambda \in I_k^n$ if and only if for all $t\in [-1,1]$,
\begin{equation}\label{eq1phi}
1 + \lambda \left(P(t) -  t P'(t)\right)\geq  0
\end{equation}
and
\begin{equation}\label{eq2phi}
1 + \lambda \left((1-t^2)P''(t) + P(t) -  t P'(t)\right )\geq  0.
\end{equation}
Now note that, by (\ref{derivLeg}), $P'(1) = k(k+n-2)/(n-1)$ and, therefore,
\begin{equation*}
P(1)-P'(1) = -\frac{(k-1)(k+n-1)}{n-1}.
\end{equation*}
Hence, by (\ref{eq1phi}), if $\lambda \in I_k$ we must have
\begin{equation*}
1 - \lambda \frac{(k-1)(k+n-1)}{n-1}\geq 0.
\end{equation*}
Rearranging we obtain the upper bound for $I_k^n$.

On the other hand, let $t_0$ be the minimizer of $P$ in $[-1,1]$. Since $t_0$ is a critical point, the derivative of $P$ at $t_0$ vanishes. Hence, evaluating (\ref{diffequpkn}) at $t_0$ yields
\begin{equation}\label{phi2}
(1-t_0^2)P''(t_0) = - k(k+n-2)P(t_0).
\end{equation}
Combining (\ref{eq2phi}) and (\ref{phi2}), yields the lower bound for $I_k^n$. Explicit computation when $k = 2$ using (\ref{p2nexpl}) shows that the bounds are optimal in this case.

Now we proceed to compute the interval $J_k^n$. Here, we just need to check for which vales of $\gamma$ the condition on Minkowski's theorem are satisfied. Note that the measure defined by the density
$1 + \gamma P(u\cdot \bar{e})$ is centered at zero by orthogonality of Legendre polynomials of different degrees. Thus we just need to check when $ 1 + \gamma P(\bar{e} \cdot u)\geq 0$ for all $u \in \mathbb{S}^{n-1}$. This happens if and only if
$$\gamma\geq \frac{-1}{\max_{[-1,1]}P_k^n} = 1 \qquad \mbox{and} \qquad \gamma \leq \frac{-1}{\min_{[-1,1]} P_k^n} = \frac{1}{\nu^n_k[1]}.$$

\vspace{-0.5cm}

\hfill $\blacksquare$

\vspace{0.4cm}

Finally, we are in a position to state and prove our new spectral gap theorem which is not only critical to establish Theorem \ref{mainthm} but might also be of independent interest in convex geometry and valuation theory.

\begin{theorem} \label{spectralgap} Suppose that $L \in \mathcal{K}^n$ is origin-symmetric and $\mathrm{SO}(n - 1)$ invariant. Then
\begin{equation*}
|a_k^n[L]| < \frac{a_0^n[L]}{(k-1)(n+k-1)}
\end{equation*}
for every $k > 2$ and
\begin{equation*}
|a_2^n[L]| \leq \frac{a_0^n[L]}{n+1},
\end{equation*}
where this inequality is also strict if $L$ is of class $C^2_+$.
\end{theorem}

\noindent {\it Proof.} Clearly, we may assume that $L$ is not a singleton and, hence, that $a_0^n[L] \neq 0$. We know that $S_{n-1}(K,\cdot) \ast h(L,\cdot)$ defines a support function for every $K \in \mathcal{K}^n$. We now choose $K = K_{\gamma}$ such that its surface area measure has a density of the form $s_{\gamma}(u) = 1 + \gamma P_k^n(\bar{e} \cdot u)$, $u \in \mathbb{S}^{n-1}$, where $\gamma$ is in the interval $J_k^n$ computed in Proposition \ref{intervals1742}. Then it follows that
\[S_{n-1}(K_\gamma,\cdot) \ast h(L,\cdot) = a_0^n[L] + a_k^n[L]\gamma P_k^n(\bar{e}\cdot .\,)   \]

\pagebreak

\noindent is the support function of a convex body in $\mathcal{K}^n$ and, hence, by Proposition \ref{intervals1742}, we must have $\gamma a_k^n[L]/a_0^n[L] \in I_k^n$ for every $\gamma \in J_k^n$ or, equivalently,
\begin{equation} \label{contained1742}
\frac{a_k^n[L]}{a_0^n[L]} J_k^n \subseteq I_k^n.
\end{equation}
Hence, if $a_k^n[L]$ is positive, then, by (\ref{decrmin}),
\[\frac{a_k^n[L]}{a_0^n[L]} \leq \frac{(n-1)\nu_k^n[1]}{(k-1)(k+n-1)} \leq \frac{1}{(k-1)(k+n-1)}.\]
for every $k \geq 2$ with strict inequality in the right hand inequality for $k > 2$. On the other hand, if $a_k^n[L]$ is negative, then
\begin{align*}
\frac{a_k^n[L]}{a_0^n[L]\nu_k^n[1]} &\geq - \frac{1}{(k(k+n-2)-1)\nu_k^n[1]}
\end{align*}
or, equivalently,
\[\frac{|a_k^n[L]|}{a_0^n[L]} \leq \frac{1}{k(k+n-2)-1}.  \]
Since $k(k+n-2)-1 > (k-1)(k+n-1)$, this proves the inequalities of the theorem. If $L$ is of class $C^2_+$, the same arguments yield the claim for $k = 2$, if we can show that $\mathrm{T}_L$ maps surface area measure densities of convex bodies to support functions of $C_2^+$ convex bodies, since in that case $I_2^n$ in (\ref{contained1742}) can be replaced by $\mathrm{int}\,I_2^n$ and both $I_2^n$ and $J_2^n$ are known precisely.

It remains to show that when $D^2(\mathrm{T}_Ls_{n-1}(K,\cdot))(u)$ is considered as a linear map on $u^{\bot}$, then $\det\left(D^2(\mathrm{T}_Ls_{n-1}(K,\cdot))(u)\right) > 0$ for all $u \in \mathbb{S}^{n-1}$ and $K \in \mathcal{K}^n$ with continuous surface area measure density. Let $\Lambda_L = \inf_{u,v\in \mathbb{S}^{n-1}} u^{\bot} D^2h(L,\cdot)(v) u>0$ denote the smallest (non-zero) eigenvalue of $D^2h(L,\cdot)$ on $\mathbb{S}^{n-1}$.
Then it follows from Example \ref{exps2} (b) and the log-concavity of the determinant that
\begin{align*}
\det\left(D^2(\mathrm{T}_Ls_{n-1}(K,\cdot))(u)\right)^{\frac{1}{n-1}} & = \det\left( \int_{\mathbb{S}^{n-1}}\!\!  D^2h(L(v),\cdot)(u) s_{n-1}(K,v)\,dv \right)^{\frac{1}{n-1}} \\
& \geq \int_{\mathbb{S}^{n-1}}\!\! \det D^2h(L(v),\cdot)(u)^{\frac{1}{n-1}} s_{n-1}(K,v)\,dv \\
& \geq \Lambda_L S(K) > 0,
\end{align*}
where $S(K)=S_{n-1}(K,\mathbb{S}^{n-1})$ denotes the surface area of $K$. \hfill $\blacksquare$

\vspace{1cm}

\centerline{\large{\bf{ \setcounter{abschnitt}{6}
\arabic{abschnitt}. Proof of the main result}}}

\reseteqn \alpheqn \setcounter{theorem}{0}

\vspace{0.6cm}

We are finally in a position to prove our main result. The ideas and techniques of Ivaki \textbf{\cite{Ivaki2018}} are the basis for the proof of the following more general, but also more technical version of Theorem \ref{mainthm}.

\begin{theorem} \label{mostgeneralthm} Let $2 \leq i \leq n - 1$ and $\Phi_i: \mathcal{K}^n \rightarrow \mathcal{K}^n$ be a continuous
translation invariant even Minkowski valuation of degree $i$ which commutes with $\mathrm{SO}(n)$ and suppose that its generating function $f \in L^1(\mathbb{S}^{n-1})$
satisfies the following conditions:
\begin{itemize}
\item[(1)] The convolution transform $\mathrm{T}_f$ is a bounded map from $C(\mathbb{S}^{n-1})$ to $C^2(\mathbb{S}^{n-1})$,
\item[(2)] for all $k \geq 2$,
\[|a_k^n[f]| < \frac{a_0^n[f]}{(k-1)(k+n-1)}\frac{n-1}{i},  \]
\item[(3)] for some integer $\varrho > 2$,
\[a_k^n[f] = \mathrm{O}(k^{-\varrho}) \qquad \mbox{ as } k \rightarrow \infty.  \]
\end{itemize}
Then there exists $\varepsilon > 0$ such that if $K \in \mathcal{K}^n$ has a $C^2$ support function and satisfies
\begin{itemize}
\item[(i)] $\|h(\gamma K + x,\cdot) - h(\mathbb{B}^n,\cdot)\|_{C^2(\mathbb{S}^{n-1})} < \varepsilon$ for some $\gamma > 0$ and $x \in \mathbb{R}^n$,
\item[(ii)] $\Phi_i^2 K = \alpha K$ for some $\alpha > 0$,
\end{itemize}
then $K$ must be a Euclidean ball.
\end{theorem}

\noindent {\it Proof.} Without loss of generality, we may assume that $\Phi_i$ is non-trivial and (as the first part of the proof of Proposition 1 shows) that we may normalize $\Phi_i$ such that $\Phi_i \mathbb{B}^n = \mathbb{B}^n$. Next, note that if $K$ satisfies (ii), then it must satisfy
\begin{equation}\label{eqfixedpointm}
        \Phi_i^{2m} K = \beta_m K
\end{equation}
for suitable $\beta_m>0$ and all integers $m \geq 1$. Taking the mean width on both sides of $(\ref{eqfixedpointm})$ yields
\begin{equation}
                \beta_m  = \frac{w(\Phi_i^{2m}K)}{w(K)}.
\end{equation}
Hence, using the projection $\pi_0$ to $\mathcal{H}_0^n$, (\ref{eqfixedpointm}) becomes in terms of support functions,
\begin{equation}
        h(\Phi_i^{2m} K,\cdot) = \frac{\int_{\mathbb{S}^{n-1}}h(\Phi_i^{2m}K,u)\,du}{\int_{\mathbb{S}^{n-1}}h(K,u)\,du}\,  h(K,\cdot) = \frac{\pi_0h(\Phi_i^{2m}K,\cdot)}{\pi_0h(K,\cdot)}\,  h(K,\cdot)
\end{equation}
for all $m \geq 1$. In the following, let $\mathcal{S}^{n,2}_{+,e}$ denote the subset of even functions in the open cone $\mathcal{S}^{n,2}_{+} \subseteq C^2(\mathbb{S}^{n-1})$. For $m \geq 1$, we define
$\mathcal{F}_m: \mathcal{S}^{n,2}_{+,e} \rightarrow C^2_e(\mathbb{S}^{n-1})$ by
\begin{equation}\label{eqF}
\mathcal{F}_m(h) = \Phi_i^{2m}(h) - \frac{\pi_0 \Phi_i^{2m}(h)}{\pi_0h} h.
\end{equation}
Here, $C^2_e(\mathbb{S}^{n-1})$ is the subspace of even function in $C^2(\mathbb{S}^{n-1})$. Clearly, $\mathcal{F}_m$ maps constants (that is, support functions of origin-symmetric balls) to zero. Our goal is to show that in a neighborhood of $h_{\mathbb{B}^n}$, constants are the only zeros of $\mathcal{F}_m$ for some $m \geq 1$. Equivalently, we need to show that $\mathcal{G}_m: \mathcal{S}^{n,2}_{+,e} \rightarrow C^2_e(\mathbb{S}^{n-1})$, defined by
\begin{equation} \label{defGm}
 \mathcal{G}_m(h)  = \mathcal{F}_m(h)   + \pi_0h,
\end{equation}
has only constant functions as fixed points.

\pagebreak

Since any zero $h$ of $\mathcal{F}_m$ is mapped to a constant by $\mathcal{G}_m$, $\mathcal{G}_m(h) = \pi_0 h = \mathcal{G}_m(\pi_0 h)$,
it will suffice to show that $\mathcal{G}_m$ is a local diffeormorphism around $h_{\mathbb{B}^n}$. In order to apply the inverse function theorem, Theorem \ref{invfctthm}, for this purpose, we first compute the Fr\'echet derivative of $\mathcal{F}_m$ at $h_{\mathbb{B}^n}$. Using basic properties of the Fr\'echet derivative, we obtain
\[\mathrm{d}\mathcal{F}_m(h)g = \mathrm{d}\Phi_i^{2m}(h)g - \frac{\pi_0\Phi_i^{2m}(h)}{\pi_0 h} g - \left(\frac{\mathrm{d}(\pi_0\Phi_i^{2m})(h)g}{\pi_0h}-\frac{\pi_0 \Phi_i^{2m}(h)\pi_0g}{{(\pi_0h)}^2}\right)h\]
for all $h\in \mathcal{S}^{n,2}_{+,e}$ and $g\in C^2_e(\mathbb{S}^{n-1})$. Since $\Phi_i^{2m}(h_{\mathbb{B}^n})= h_{\mathbb{B}^n}$ and $\pi_0h_{\mathbb{B}^n} = 1$, it follows that
\begin{align}\label{diff1}
\mathrm{d}\mathcal{F}_m(h_{\mathbb{B}^n})g = \mathrm{d}\Phi_i^{2m}(h_{\mathbb{B}^n})g- g - \mathrm{d}(\pi_0\Phi_i^{2m})(h_{\mathbb{B}^n})g + \pi_0g.
\end{align}
Since $\pi_0$ is linear and bounded, the chain rule, followed by an application of Corollary~\ref{diffval}, yields
\[\mathrm{d}(\pi_0 \Phi_i^{2m})(h_{\mathbb{B}^n})g = \pi_0 \mathrm{d}\Phi_i^{2m}(h_{\mathbb{B}^n})g = \pi_0(i\,\Box_n \mathrm{T}_f)^{2m}g = i^{2m} \pi_0g,\]
where for the last equality, we have used that $a_0^n[f] = 1$, by our normalization of $\Phi_i$, and $a_0^n[\Box_n] = 1$, by (\ref{boxnmult}).
Substituting in (\ref{diff1}) and using again Corollary \ref{diffval}, we obtain
\begin{equation} \label{useful1742}
\mathrm{d}\mathcal{F}_m(h_{\mathbb{B}^n})g  =  (i\,\Box_n\mathrm{T}_f)^{2m} g - g - (i^{2m} - 1)\pi_0g
\end{equation}
for all $g\in C^2_e(\mathbb{S}^{n-1})$. Next, we want to use (\ref{useful1742}) to determine the kernel of $\mathrm{d}\mathcal{F}_m(h_{\mathbb{B}^n})$. Indeed, we claim that
\begin{equation}
\ker \mathrm{d}\mathcal{F}_m(h_{\mathbb{B}^n}) = \mathcal{H}_0^n.
\end{equation}
To see this, first note that, by (\ref{useful1742}) and the fact that $a_0^n[f] = 1$ and $a_0^n[\Box_n] = 1$, we have on one hand that $\pi_0 \mathrm{d}\mathcal{F}_m(h_{\mathbb{B}^n})g = 0$ for all $g\in C^2_e(S^{n-1})$.
On the other hand, by (\ref{useful1742}) and (\ref{boxnmult}), we have for every $k\geq 1$,
\begin{equation}
\pi_{2k} \mathrm{d}\mathcal{F}_m(h_{\mathbb{B}^n})g = \left (i^{2m}\frac{(1-2k)^{2m}(2k+n-1)^{2m}}{(n-1)^{2m}}a_{2k}^n[f]^{2m} - 1\right)\pi_{2k} g
\end{equation}
for all $g\in C^2_e(\mathrm{S}^{n-1})$. Applying Parseval's identity, we obtain
\begin{equation*}
 \|\mathrm{d}\mathcal{F}_m(h_{\mathbb{B}^n})g\|_{L^2}^2 = \sum_{k = 1}^\infty \left(i^{2m}\frac{(1-2k)^{2m}(2k+n-1)^{2m}}{(n-1)^{2m}}a_{2k}^n[f]^{2m} - 1\right)^{\!2}\|\pi_{2k} g\|_{L^2}^2
\end{equation*}
for all $g\in C^2_e(\mathbb{S}^{n-1})$. By assumption (2) on the multipliers $a_{k}^n[f]$, all the leading coefficients of this sum are non-zero. Thus, $g \in \ker{\mathrm{d}\mathcal{F}_m(h_{\mathbb{B}^n})}$ if and only if $\pi_{2k}g = 0$ for all $k\geq 1$. Consequently, $\ker{\mathrm{d}\mathcal{F}_m(h_{\mathbb{B}^n})} = \mathcal{H}^n_0$.
In particular, it follows that $\mathrm{d}\mathcal{G}_m$ is injective in a neighborhood of $h_{\mathbb{B}^n}$.

In order to prove surjectivity of $\mathrm{d}\mathcal{G}_m$ at $h_{\mathbb{B}^n}$, it suffices, by (\ref{defGm}) and (\ref{useful1742}), to show that for each $h\in C^2_e(\mathbb{S}^{n-1})$ such that $\pi_0h = 0$, there exists a unique $g\in C^2_e(\mathbb{S}^{n-1})$ with
$\pi_0g = 0$ such that
\begin{equation}\label{solveq}
(i\,\Box_n \mathrm{T}_{f})^{2m} g - g = h.
\end{equation}
By assumption (3), we have
\[\lim_{k\to\infty}(1-2k)^{2m}(2k+n-1)^{2m}a_{2k}^n[f]^{2m} = 0.\]
Combining this with (2), we conclude that the series
\[\sum_{k = 0}^\infty \left (i^{2m}\frac{(1-2k)^{2m}(2k+n-1)^{2m}}{(n-1)^{2m}}a_{2k}^n[f]^{2m} - 1   \right )^{\!-1}\pi_k h\]
convergences to an even function $\zeta \in L^2(\mathbb{S}^{n-1})$. Define $g =  (i\,\Box_n\mathrm{T}_{f})^{2m} \zeta - h$.
From (\ref{smoothm}) and the conclusion following it together with (3) and the fact that $a_k^n[\Box_n] = \mathrm{O}(k^2)$, we deduce that  $(i\,\Box_n\mathrm{T}_f)^{2m} \zeta \in H_e^{2m(\varrho-2)}(\mathbb{S}^{n-1})$. Thus, by (\ref{sobembed}),
$$
(i\,\Box_n\mathrm{T}_f)^{2m} \zeta \in C^2_e(\mathbb{S}^{n-1})
$$
provided that $m > \frac{n+3}{4(\varrho-2)}$ and, hence, $g \in C^2_e(\mathbb{S}^{n-1})$. Finally, by construction, $g$ satisfies the desired equation (\ref{solveq}).

By Theorem \ref{invfctthm} applied to the map $\mathcal{G}_m$, there exists a $C^2$ neighborhood $B_\varepsilon$ of $h_{\mathbb{B}^n}$, where $\mathcal{G}_{m}$ is a diffeomorphism.
Thus, if $K \in \mathcal{K}^n$ has $C^2$ support function and satisfies (i) and (ii), then
\[\mathcal{G}_m(h(K,\cdot)) = \mathcal{F}_m(h(K,\cdot)) + \pi_0h(K,\cdot)  = \pi_0h(K,\cdot) = \mathcal{G}_m(\pi_0h(K,\cdot)).\]
Moreover, since
\begin{equation*}
  |h(\mathbb{B}^n,\cdot) - \pi_0h(K,\cdot)| \leq \frac{1}{\omega_{n-1}}\int_{\mathbb{S}^{n-1}} \!\! |1 - h(K,u)|\, du \leq \|h(\mathbb{B}^n,\cdot)-h(K,\cdot)\|_{C^2(\mathbb{S}^{n-1})}.
\end{equation*}
Hence,  $\pi_0h(K,\cdot) \in B_\varepsilon$. Since $\mathcal{G}_m$ is bijective on $B_\varepsilon$, we have that $h(K,\cdot) = \pi_0 h(K,\cdot)$ and so $K$ is a ball. \hfill $\blacksquare$

\vspace{0.3cm}

Before we show how Theorem \ref{mainthm} can be deduced from Theorem \ref{mostgeneralthm}, we want to make a couple of remarks about the above proof.
First, note that the assumption that $\Phi_i$ is even was crucial since there is no analogue of Strichartz' result (\ref{strichartz}) for general functions on $\mathbb{S}^{n-1}$ (as far as the authors are aware).
Second, the arguments in the proof of Theorem \ref{mostgeneralthm} also hold when $\Phi_i^2$ in the fixed point assumption (ii) is replaced by $\Phi_i$. Finally, we note that the additional iterations $\Phi_i^{2m}$ for $m \geq 1$ are not required if $f$ is sufficiently regular, for example, when $f$ is smooth.

\vspace{0.2cm}

Let us show now how Theorem \ref{mainthm} follows from Theorem \ref{mostgeneralthm} using Theorem \ref{spectralgap} almost effortlessly.

\vspace{0.3cm}

\noindent {\it Proof of Theorem \ref{mainthm}.} By Theorem \ref{mostgeneralthm}, it suffices to show that if $L \in \mathcal{K}^n$ is an origin-symmetric body of revolution of class $C^2_+$, then the convolution transform $\mathrm{T}_L$ satisfies the conditions (1), (2), and (3) from Theorem \ref{mostgeneralthm}. However, (1) was already proved in Example \ref{exps2} (b) and (2) is the content of Theorem \ref{spectralgap}.

In order to show (3), note that, by (\ref{multconvpkn}), (\ref{deltasmult}), and the fact that the spherical Laplacian $\Delta_{\mathbb{S}}$ is self-adjoint,
\begin{equation*}
a_k^n[L]  = -\frac{1}{k(k+n-2)} \int_{\mathbb{S}^{n-1}}\!\! \Delta_{\mathbb{S}}h(L,\cdot)(u)P_k^n(\bar{e}\cdot u)\,du.
\end{equation*}
Hence, by the Cauchy--Schwarz inequality,
 \begin{equation*}
|a_k^n[L]|  \leq \frac{1}{k^2}\sqrt{\frac{\omega_n}{N(n,k)}}{||\Delta_{\mathbb{S}} h(L,\cdot)||}_{L^2} < \infty.		
\end{equation*}
Consequently, by (\ref{nnk}), we obtain the desired asymptotic estimate
\begin{equation*}
a_k^n[L] =  \mathrm{O}\!\left(k^{-\frac{n+2}{2}}\right) \quad \mbox{as } k \to \infty.
\end{equation*}

\vspace{-0.5cm}

\hfill $\blacksquare$

\vspace{0.4cm}

We state one more consequence of Theorem \ref{mostgeneralthm} which shows that $C^2_+$ regularity can be relaxed when $\Phi_i$ is generated by a generalized zonoid of revolution.

\begin{koro} \label{genzoncor} Let $n \geq 4$, $2 \leq i \leq n - 2$, and $\Phi_i: \mathcal{K}^n \rightarrow \mathcal{K}^n$ be a continuous translation invariant even Minkowski valuation of degree $i$ which commutes with $\mathrm{SO}(n)$. If $\Phi_i$ is generated by generalized zonoid of revolution $Z^\mu \in \mathcal{K}^n$, then there exists $\varepsilon > 0$ such that if $K \in \mathcal{K}^n$ has a $C^2$ support function and satisfies
\begin{itemize}
\item[(i)] $\|h(\gamma K + x,\cdot) - h(\mathbb{B}^n,\cdot)\|_{C^2(\mathbb{S}^{n-1})} < \varepsilon$ for some $\gamma > 0$ and $x \in \mathbb{R}^n$,
\item[(ii)] $\Phi_i^2 K = \alpha K$ for some $\alpha > 0$,
\end{itemize}
then $K$ must be a Euclidean ball.

\end{koro}

\noindent {\it Proof.} By Theorem \ref{mostgeneralthm}, it suffices to show that the convolution transform $\mathrm{T}_{Z^{\mu}}$ satisfies the conditions (1), (2), and (3) from Theorem \ref{mostgeneralthm}. Condition (2) is again a consequence of Theorem \ref{spectralgap} since $i \leq n - 2$. Thus, it remains to show that conditions (1) and (3) are satisfied.
To this end, first note that, by (\ref{TZmu}), $\mathrm{T}_{Z^{\mu}} = \mathrm{C}\circ \mathrm{T}_{\mu}$.

In order to show that $\mathrm{T}_{Z^\mu}$ satisfies (1), we identify the even signed zonal measure $\mu$ on $\mathbb{S}^{n-1}$ with an even signed measure $\bar{\mu}$ on $[-1,1]$ and use the fact (see, e.g., \textbf{\cite{kiderlen05}}) that for any $f\in L^2(\mathbb{S}^{n-1})$, $\mathrm{T}_{\mu}f$ can be written in the form
\[ \mathrm{T}_\mu f(u) =  \int_{\!-1}^1 \mathrm{R}_{\,t} f(u)\, d\bar{\mu}(t), \qquad u \in \mathbb{S}^{n-1},\]
where $\mathrm{R}_{\,t}$ denotes the generalized spherical Radon transforms defined by
\[ \mathrm{R}_{\,t}f(u) = \frac{1}{\omega_{n-2}}\int_{\mathbb{S}^{n-1}\cap u^\bot}\!\!\! f\left(t u + \sqrt{1-t^2}v\right)\,dv.\]
Thus, $\mathrm{T}_{Z^{\mu}}f = \int_{\!-1}^1 \mathrm{C}\mathrm{R}_{\,t} f\, d\bar{\mu}(t)$.
Since $\bar{\mu}$ is even, $\mathrm{R}_{\,1} = \mathrm{Id}$, $\mathrm{R}_{\,-1} = \mathrm{-Id}$, and we know from Example \ref{exps2} (a) that $\mathrm{C}$ satisfies condition (1), it suffices, by the uniform boundedness principle, to prove that $\mathrm{C}\mathrm{R}_{\,t}$ satisfies (1) for all $t \in (-1,1)$. To see this, we use the fact (see, e.g., \textbf{\cite{kiderlen05}}) that for all $t \in [-1,1]$,
\[a_k^n[\mathrm{R}_{\,t}] = P_k^n(t),\]
combined with the classical asymptotic estimate (see, e.g., \textbf{\cite[\textnormal{p.\ 172}]{Szegoe1975}})
\[P_k^n(t) = \arccos(t)^{-\frac{n-2}{2}}\, \mathrm{O}\!\left(k^{-\frac{n-2}{2}}\right) \quad \mbox{as } k \to \infty\]
for all $t \in (-1,1)$. From this and the multipliers of the cosine transform (\ref{multC}), it follows that $a_k^n[\mathrm{C}\mathrm{R}_{\,t}]=a_k^n[\mathrm{C}]a_k^n[\mathrm{R}_{\,t}]=  \mathrm{O}(k^{-n})$ as $k \to \infty$. Hence, by the smoothing property (\ref{smoothm}) and the Sobolev embedding (\ref{sobembed}), the operators $\mathrm{C}\mathrm{R}_{\,t}$ are bounded from $L^2(\mathbb{S}^{n-1}) \to C^2(\mathbb{S}^{n-1})$ for all $t \in (-1,1)$ as long as $n\geq 4$.

In order to prove that $\mathrm{T}_{Z^\mu}$ satisfies (3), note that $a_k^n[\mathrm{T}_{Z^\mu}]=a_k^n[\mathrm{C}]a_k^n[\mu]$. Therefore, by (\ref{multconvpkn}) and the fact that $|P_k^n| \leq 1$,
\[|a_k^n[\mathrm{T}_{Z^\mu}]| = ||\mu||_{\mathrm{TV}}  \mathrm{O}\!\left(k^{-\frac{n+2}{2}}\right) \quad \mbox{as } k \to \infty,\]
where  $||\mu||_{\mathrm{TV}} = \mu^+(\mathbb{S}^{n-1}) + \mu^{-}(\mathbb{S}^{n-1})<\infty$.
\hfill $\blacksquare$

\vspace{0.3cm}

We conclude the article with the remark that if $i = n - 1$ and $Z^{\mu}$ is generated by a \emph{non-negative} and \emph{non-discrete} measure $\mu$, then the statement of Corollary \ref{genzoncor} holds true by (essentially) the same arguments as given above.

\vspace{0.5cm}

\noindent {{\bf Acknowledgments} The authors were supported by the Austrian Science Fund (FWF), Project number: P31448-N35.

\begin{small}

\[ \begin{array}{ll} \mbox{Oscar Ortega Moreno} & \mbox{Franz E. Schuster} \\
\mbox{Vienna University of Technology \phantom{wwwwWW}} & \mbox{Vienna University of Technology} \\ \mbox{oscarortem@gmail.com} & \mbox{franz.schuster@tuwien.ac.at}
\end{array}\]

\end{small}

\end{document}